\documentclass[onefignum,onetabnum]{siamart190516}

\usepackage{algpseudocode}
\usepackage{amsmath,amsfonts}
\usepackage{bbm}
\usepackage[sort]{cite}
\usepackage{cleveref}
\usepackage{color}
\usepackage{csl}
\usepackage{dsfont}
\usepackage{forest}
\usepackage{graphicx}
\usepackage{mathtools}
\usepackage{mleftright}
\usepackage{physics}
\usepackage{pgfplots}
\usepackage{pgfplotstable}
\usepackage{siunitx}
\usepackage{subcaption}
\usepackage{xparse}

\graphicspath{{./Figures/}}
\newcommand{\datadir}{./Data}

\pgfplotsset{compat=1.14}
\pgfplotscreateplotcyclelist{mycolor1}{%
	red,mark=*\\%
	blue,mark=triangle\\%
	black,thick,mark=star\\%
	black,densely dashed,mark=star\\%
}
\pgfplotscreateplotcyclelist{mycolor2}{%
	red,mark=*\\%
	red,densely dashed,mark=*\\%
	blue,mark=triangle\\%
	blue,densely dashed,mark=triangle\\%
	black,thick,mark=star\\%
	black,densely dashed,mark=star\\%
}
\pgfplotsset{every axis/.append style={
	label style={font=\scriptsize},
	legend style={font=\scriptsize,/tikz/every even column/.append style={column sep=1em}},
	tick label style={font=\scriptsize},
	legend cell align={left},
	yminorticks=false,
	grid,
	width=\linewidth
}}
\usepgfplotslibrary{groupplots}

\newsiamremark{remark}{Remark}

\NewDocumentCommand{\mrimethod}{}{SM-MRI-GARK}
\NewDocumentCommand{\spcmethod}{}{SM-SPC-MRI-GARK}

\NewDocumentCommand{\eye}{g g}{\mathbf{I}\IfValueT{#1}{_{{#1} \times \IfValueTF{#2}{#2}{#1}}}}
\NewDocumentCommand{\zero}{g g}{\IfValueTF{#2}{\mathbf{0}}{0}\IfValueT{#1}{_{{#1} \IfValueT{#2}{\times {#2}}}}}
\NewDocumentCommand{\one}{g g}{\mathbbm{1}\IfValueT{#1}{_{{#1} \IfValueT{#2}{\times {#2}}}}}

\NewDocumentCommand{\R}{g g g}{\mathbb{R}^{\IfValueT{#1}{#1} \IfValueT{#2}{\times {#2}} \IfValueT{#3}{\times {#3}}}}
\NewDocumentCommand{\Cplx}{g g g}{\mathbb{C}^{\IfValueT{#1}{#1} \IfValueT{#2}{\times {#2}} \IfValueT{#3}{\times {#3}}}}

\DeclareMathAlphabet{\mathbase}{OT1}{pzc}{m}{n}
\NewDocumentCommand{\Abase}{O{} s}{\mathbase{A}\component{\slow,\slow}{#1}{#2}}
\NewDocumentCommand{\abase}{O{} s}{\mathbase{a}\component{\slow}{#1}{#2}}
\NewDocumentCommand{\bbase}{O{} s}{\mathbase{b}\component{\slow}{#1}{#2}}
\NewDocumentCommand{\bhatbase}{O{} s}{\widehat{\mathbase{b}}\component{\slow}{#1}{#2}}
\NewDocumentCommand{\cbase}{O{} s}{\mathbase{c}\component{\slow}{#1}{#2}}
\NewDocumentCommand{\deltac}{O{} s}{{\Delta \mathbase{c}}\component{\slow}{#1}{#2}}
\NewDocumentCommand{\scalarstab}{O{\rho} O{\alpha}}{\mathcal{S}^{\textsc{1d}}_{#1,#2}}

\newcommand{\fast}{\mathfrak{f}}
\newcommand{\slow}{\mathfrak{s}}
\NewDocumentCommand{\component}{m m m}{^{ \left\{ #1 \right\} #2 \IfBooleanT{#3}{T} }}
\NewDocumentCommand{\comp}{m O{} s}{\component{#1}{#2}{#3}}
\NewDocumentCommand{\F}{O{} s}{\component{\fast}{#1}{#2}}
\NewDocumentCommand{\FL}{O{\lambda} O{} s}{\component{\fast, #1}{#2}{#3}}
\RenewDocumentCommand{\S}{O{} s}{\component{\slow}{#1}{#2}}
\NewDocumentCommand{\SL}{O{\lambda} O{} s}{\component{\slow, #1}{#2}{#3}}
\NewDocumentCommand{\FF}{O{} s}{\component{\fast, \fast}{#1}{#2}}
\NewDocumentCommand{\FS}{O{} s}{\component{\fast, \slow}{#1}{#2}}
\NewDocumentCommand{\FSL}{O{\lambda} O{} s}{\component{\fast, \slow, #1}{#2}{#3}}
\NewDocumentCommand{\SF}{O{} s}{\component{\slow, \fast}{#1}{#2}}
\NewDocumentCommand{\SFL}{O{\lambda} O{} s}{\component{\slow, \fast, #1}{#2}{#3}}
\RenewDocumentCommand{\SS}{O{} s}{\component{\slow, \slow}{#1}{#2}}

\NewDocumentCommand{\qgconvplot}{m m s}{
	\begin{tikzpicture}
	\begin{loglogaxis}[
		cycle list name=mycolor2,
		height=0.4\linewidth,
		legend columns=2,
		legend entries={
			\mrimethod{} using ROM of size 80,
			\mrimethod{} using ROM of size 40,
			\spcmethod{} using ROM of size 80,
			\spcmethod{} using ROM of size 40,
			Runge--Kutta using full model},
		\IfBooleanTF{#3}{xlabel=Steps,legend to name=legend:qg_conv}{legend entries={}},
		ylabel={Relative $\ell^2$ Error}]
		
		\pgfplotstableread[col sep=comma]{\datadir/QG/exp_qgrom_gark_#1.dat}{\qggark}
		\pgfplotstableread[col sep=comma]{\datadir/QG/exp_qgrom_spc_#1.dat}{\qgspc}
		\addplot table[x=Steps,y=romr80] {\qggark};
		\addplot table[x=Steps,y=romr40] {\qggark};
		\addplot table[x=Steps,y=romr80] {\qgspc};
		\addplot table[x=Steps,y=romr40] {\qgspc};
		\addplot table[x=Steps,y=reference] {\qgspc};
	
		\pgfplotstablecreatecol[linear regression]{regression}{\qggark}
		\pgfmathsetmacro{\slope}{-\pgfplotstableregressiona}
		\pgfplotstablegetelem{2}{Steps}\of{\qggark}
		\pgfmathsetmacro{\xleft}{\pgfplotsretval}
		\pgfplotstablegetelem{4}{Steps}\of{\qggark}
		\pgfmathsetmacro{\xright}{\pgfplotsretval}
		\draw (axis cs: \xleft,#2) -- node[below]{\scriptsize{\pgfmathprintnumber[precision=0]{\slope}}} (axis cs: \xright,{#2 * (\xleft / \xright)^\slope});
	\end{loglogaxis}
	\end{tikzpicture}}

\NewDocumentCommand{\mlconvplot}{m m s}{
	\begin{tikzpicture}
	\begin{loglogaxis}[
		cycle list name=mycolor1,
		legend columns=1,
		legend entries={\mrimethod, \spcmethod, Runge--Kutta using full model},
		\IfBooleanTF{#3}{legend to name=legend:ml_conv}{legend entries={}},
		xlabel=Steps,
		ylabel={Relative $\ell^2$ Error}]
		
		\pgfplotstableread[col sep=comma]{\datadir/L96ML/exp_mlrom_gark_#1.dat}{\mlgark}
		\pgfplotstableread[col sep=comma]{\datadir/L96ML/exp_mlrom_spc_#1.dat}{\mlspc}
		\addplot table[x=Steps,y=romr] {\mlgark};
		\addplot table[x=Steps,y=romr] {\mlspc};
		\addplot table[x=Steps,y=reference]	{\mlgark};
		
		\pgfplotstablecreatecol[linear regression={y=romr}]{regression}{\mlgark}
		\pgfmathsetmacro{\slope}{-\pgfplotstableregressiona}
		\pgfplotstablegetelem{2}{Steps}\of{\mlgark}
		\pgfmathsetmacro{\xleft}{\pgfplotsretval}
		\pgfplotstablegetelem{5}{Steps}\of{\mlgark}
		\pgfmathsetmacro{\xright}{\pgfplotsretval}
		
		\draw (axis cs: \xleft,#2) -- node[above]{\scriptsize{\pgfmathprintnumber[precision=0]{\slope}}} (axis cs: \xright,{#2 * (\xleft / \xright)^\slope});
	\end{loglogaxis}
	\end{tikzpicture}}

\NewDocumentCommand{\mladaptplot}{m}{
	\begin{tikzpicture}
	\begin{axis}[
	xlabel=Step,
	ylabel={$H$}]
	
	\addplot[
	scatter,
	scatter src=explicit,
	scatter/classes={0={mark=*,red},1={mark=none}}] table[
	col sep=comma,
	x=Step,
	y=H,
	meta=Accepted] {\datadir/L96ML_Adaptive/exp_mlrom_#1.dat};
	\end{axis}
	\end{tikzpicture}}

\NewDocumentCommand{\mltolplot}{}{
	\begin{tikzpicture}
		\begin{loglogaxis}[
			cycle list name=mycolor2,
			height=0.4\linewidth,
			xlabel={$\texttt{AbsTol}=\texttt{RelTol}$},
			ylabel={Global Relative $\ell^2$ Error},
			legend entries={Order 2 \mrimethod{},Order 3 \mrimethod{},Order 2 \spcmethod{},Order 3 \spcmethod{},Ralston 2, Ralston 3},
			legend columns=2,
			legend to name=legend:ml_tol]
			\foreach \n in {1,...,6} {
				\addplot table[col sep=comma,x=Tolerance,y index=\n] {\datadir/L96ML_Adaptive/tolerances.dat};
			}
		\end{loglogaxis}
	\end{tikzpicture}}

\NewDocumentCommand{\brusperfplot}{m s}{%
	\begin{tikzpicture}
	\begin{loglogaxis}[
		cycle list name=mycolor2,
		legend columns=2,
		height=0.4\linewidth,
		legend entries={
			\shortstack[l]{\mrimethod{} using $129 \times 129$ surrogate\\model},
			\shortstack[l]{\mrimethod{} using $65 \times 65$ surrogate\\model},
			\shortstack[l]{\spcmethod{} using $129 \times 129$\\surrogate model},
			\shortstack[l]{\spcmethod{} using $65 \times 65$\\surrogate model},
			Traditional Runge--Kutta using full model},
		\IfBooleanTF{#2}{xlabel={Runtime (s)},legend to name=legend:brus_perf}{legend entries={}},
		ylabel={Absolute $\ell^2$ Error}]
		\pgfplotstableread[col sep=comma]{\datadir/Brusselator/order#1.dat}{\brusdata}
		\addplot table[x=MRI-GARK 129 Time,y=MRI-GARK 129 Error] {\brusdata};
		\addplot table[x=MRI-GARK 65 Time,y=MRI-GARK 65 Error] {\brusdata};
		\addplot table[x=SPC-MRI-GARK 129 Time,y=SPC-MRI-GARK 129 Error] {\brusdata};
		\addplot table[x=SPC-MRI-GARK 65 Time,y=SPC-MRI-GARK 65 Error] {\brusdata};
		\addplot table[x=Traditional 257 Time,y=Traditional 257 Error] {\brusdata};
	\end{loglogaxis}
	\end{tikzpicture}}                            

\NewDocumentCommand{\advperfplot}{m s}{%
	\begin{tikzpicture}
	\begin{loglogaxis}[
		cycle list name=mycolor1,
		legend columns=1,
		legend entries={\mrimethod,\spcmethod,Runge--Kutta using full model,Runge--Kutta using surrogate model},
		\IfBooleanTF{#2}{legend to name=legend:adv_perf}{legend entries={}},
		xlabel={Runtime (s)},
		ylabel={Absolute $L^2$ Error}]
	\pgfplotstableread[col sep=comma]{\datadir/Advection/order#1.dat}{\advdata}
	\addplot table[x=MRI-GARK Time,y=MRI-GARK Error] {\advdata};
	\addplot table[x=SPC-MRI-GARK Time,y=SPC-MRI-GARK Error] {\advdata};
	\addplot table[x=Traditional Time,y=Traditional Error] {\advdata};
	\addplot table[x=Surrogate Time,y=Surrogate Error] {\advdata};
	\end{loglogaxis}
	\end{tikzpicture}}

\forestset{
	btree/.style={
		for tree={
			grow=north,
			minimum width=5pt,
			minimum height=5pt,
			inner sep=0pt,
			outer sep=0pt,
			s sep=3pt,
			circle,
			fill
		},
		before computing xy={
			for tree={
				l=8pt,
			}
		},
		delay={
			for tree={
				if={strequal("f",content)}{fill=black}{},
				if={strequal("s",content)}{fill=white,draw}{},
				content={}
			}
		}
	}
}
\newcommand{\btree}[1]{\Forest{btree #1}}

\newcommand{\nvar}{{\rm N}}
\newcommand{\nsur}{{\rm S}}
\newcommand{\ysur}{\y_{\rm sur}}
\newcommand{\fsur}{\fun_{\rm sur}}
\NewDocumentCommand{\Vmat}{s}{\mathbf{V}\IfBooleanT{#1}{^*}}
\NewDocumentCommand{\Wmat}{s}{\mathbf{W}\IfBooleanT{#1}{^*}}

\newenvironment{smatrix}[1][1.25]{\def\arraystretch{#1}\begin{bmatrix}}{\end{bmatrix}}
\newenvironment{butchertableau}[2][1.3]{\def\arraystretch{#1}\array{#2}}{\endarray}

\newcommand{\fun}{f}
\newcommand{\Jac}{\mathbf{J}}
\newcommand{\Lb}{\mathbf{L}}
\newcommand{\y}{y}

\newcommand{\sfrac}[2]{\mbox{\footnotesize$\displaystyle\frac{#1}{#2}$}} 

\newcommand{\utree}{{\mathfrak{u}}}
\newcommand{\BT}{\mathds{T}}

\newcommand{\papertitle}{A fast time-stepping strategy for dynamical systems equipped with a surrogate model}
\headers{Time-stepping with a surrogate model}{S. Roberts, A. A Popov, A. Sarshar, and A. Sandu}
\title{\papertitle{}\thanks{Submitted to the editors DATE.
		\funding{This work was supported by awards NSF ACI--1709727, NSF CDS\&E--MSS 1953113, DE-SC0021313, and by the Computational Science Laboratory at Virginia Tech.}}}
\author{Steven Roberts\thanks{\csl{}, \cslinstitution{}, \csladdress{}
		(\email{steven94@vt.edu}, \email{apopov@vt.edu}, \email{sarshar@vt.edu}, \email{sandu@cs.vt.edu}).}
	\and Andrey A Popov\footnotemark[2]
	\and Arash Sarshar\footnotemark[2]
	\and Adrian Sandu\footnotemark[2]}

\newif\ifreport
\reporttrue

\begin{document}
	
\ifreport
	\csltitle{\papertitle{}}
	\cslauthor{Steven Roberts, Andrey A Popov, Arash Sarshar, and Adrian Sandu}
	\cslemail{steven94@vt.edu, apopov@vt.edu, sarshar@vt.edu, sandu@cs.vt.edu}
	\cslreportnumber{6}
	\cslyear{20}
	\csltitlepage
\fi

\maketitle

\begin{abstract}
	Simulation of complex dynamical systems arising in many applications is computationally challenging due to their size and complexity.  Model order reduction, machine learning, and other types of surrogate modeling techniques offer cheaper and simpler ways to describe the dynamics of these systems but are inexact and introduce additional approximation errors.  In order to overcome the computational difficulties of the full complex models, on one hand, and the limitations of surrogate models, on the other, this work proposes a new accelerated time-stepping strategy that combines information from both. This approach is based on the multirate infinitesimal general-structure additive Runge--Kutta (MRI-GARK) framework.  The inexpensive surrogate model is integrated with a small timestep to guide the solution trajectory, and the full model is treated with a large timestep to occasionally correct for the surrogate model error and ensure convergence.  We provide a theoretical error analysis, and several numerical experiments, to show that this approach can be significantly more efficient than using only the full or only the surrogate model for the integration.
\end{abstract}




\section{Introduction}

This paper will consider the system of ordinary differential equations (ODEs)
\begin{equation} \label{eq:ode}
	\y' = \fun(t, \y),
	\quad
	\y(t_0) = \y_0,
	\quad
	\y \in \Cplx{\nvar},
\end{equation}
equipped with a surrogate model approximating the dynamics of \cref{eq:ode}:
\begin{equation} \label{eq:surrogate_ode}
	\ysur' = \fsur(t, \ysur),
	\qquad
	\ysur \in \Cplx{\nsur}.
\end{equation}
This surrogate model is assumed to be much less expensive to solve than the full model \cref{eq:ode}, possibly by evolving in a lower dimensional space ($\nsur < \nvar$).  Moreover, we assume that projections between the full and surrogate model spaces are realized by the matrices $\Vmat, \Wmat \in \Cplx{\nvar}{\nsur}$:
\begin{equation} \label{eq:space_transformation}
	\ysur = \Wmat* \y,
	\quad
	\y \approx \Vmat \ysur,
	\quad
	\Wmat* \Vmat = \eye{\nsur}.
\end{equation}
Here, the ``$\ast$'' symbol denotes the conjugate transpose of a matrix.

This paper presents a new technique to supplement the numerical integration of \cref{eq:ode} using the surrogate model \cref{eq:surrogate_ode}.  For an appropriate choice of surrogate model, our method can be significantly more efficient than using either the full or the surrogate model alone.  This is accomplished by applying a multirate time integration scheme to the following ODE with fast dynamics $f\F$ and slow dynamics $f\S$:
\begin{equation} \label{eq:full_surrogate_ode}
	\y' = \underbrace{\Vmat \fsur \mleft(t, \Wmat^* \y \mright)}_{\fun\F(t, \y)} + \underbrace{\fun(t, \y) - \Vmat \fsur \mleft(t, \Wmat^* \y \mright)}_{\fun\S(t, \y)}.
\end{equation}
Multirate methods are characterized by using different stepsizes for different parts of an ODE, as opposed to a single, global timestep \cite{Rice_1960_splitRK,Gear_1974_MR,Andrus_1979_MR,Gear_1984_MR-LMM,Sandu_2007_MR-monotonic-RK2,Sandu_2016_MR-GARK,Sandu_2009_MR-monotonic-LMM,Sandu_2013_MR-extrapolation,Sandu_2019_MR-GARK_High-Order,Sandu_2021_MR-GARK-ROS,Sandu_2021_MR-GARK_Implicit,Sandu_2021_MR-Euler-DAE}.  For problems with partitions exhibiting vastly different time scales, stiffnesses, evaluation costs, or amounts of nonlinearity, multirate methods can be more efficient that their single rate counterparts.  In \cref{eq:full_surrogate_ode}, the fast partition $f\F$ is treated with a small timestep, and the slow partition $f\S$ is treated with a large timestep.  Note that \cref{eq:full_surrogate_ode} has the same solution as \cref{eq:ode}.  It is rewritten and partitioned, however, in such a way that $f\F$ contains the surrogate model dynamics and $f\S$ represents the error in the surrogate model.  Ideally, this error will be small, so a large timestep would be acceptable.  Moreover, this means expensive evaluations of the full model occur infrequently compared to the inexpensive surrogate model evaluations.

With about 60 years of development \cite{Rice_1960_splitRK}, the multirating strategy has been applied to numerous classes of time integration methods.  Conceivably any multirate method suitable for additively partitioned systems could be applied to \cref{eq:full_surrogate_ode}.  This paper builds upon multirate infinitesimal methods \cite{Knoth_1998_MR-IMEX,Schlegel_2009_RFSMR,Wensch_2009_MIS,sexton2018relaxed,Bauer2019}, as they offer a particularly flexible and elegant approach to multirate integration.  In particular, we use the multirate infinitesimal general-structure additive Runge--Kutta (MRI-GARK) framework proposed in \cite{Sandu_2019_MRI-GARK,Sandu_2020_MRI-GARK_Coupled}.  MRI-GARK is an appealing choice as it generalizes many types of multirate infinitesimal methods, allows for the construction of high order methods, and supports implicit stages.  We note, however, that this paper primarily focuses on explicit methods where \cref{eq:ode} is nonstiff.

There are other instances in the literature where multirate methods and surrogate models are used in conjunction.  In \cite{hachtel2018multirate}, an implicit multirate scheme is used to simulate a fast-evolving electric circuit coupled with a slow-evolving electric field.  Model order reduction is then applied to the electric field problem to further reduce the computational cost.  A multirate Runge--Kutta--Chebyshev method is proposed in \cite{abdulle2020stabilized} using a time-averaged right-hand side, which can be viewed as a surrogate model with better stiffness properties.  The idea of treating models of different fidelities with different timesteps is also explored in \cite{keshavarzzadeh2019convergence,liu2020hierarchical}. Surrogate models have been used extensively to speed up optimization problems when objective function evaluations are expensive \cite{ong2003evolutionary,queipo2005surrogate,forrester2009recent}.

From the multiscale modeling community, there are several related ideas such as high-order/low-order methods \cite{chacon2017multiscale}.  Similar to the methods proposed in this work, the heterogeneous multiscale method and projective integration \cite{vanden2007hmm,abdulle_weinan_engquist_vanden-eijnden_2012} utilize different timesteps for ``microscopic'' and ``macroscopic'' models as well as operators to convert a state between the two models.  They seek to efficiently approximate macroscopic quantities of interest with the help of a microscopic model.  In their implementations, the microscopic model is evaluated in small bursts of microsteps to resolve transient effects, while the macroscopic model is evolved with large timesteps.  Typically, the coupling between models introduces first order temporal errors.  In contrast, the new schemes proposed in this work seek to approximate the full, microscopic model to high orders using large macrosteps.

This remainder of this paper is structured as follows.  In \cref{sec:mri_gark}, we review the MRI-GARK framework and extensions.  \Cref{sec:method_formulation} then specializes MRI-GARK-type methods to the special ODE \cref{eq:full_surrogate_ode}.  An error analysis is performed in \cref{sec:error}, and various approaches for constructing surrogate models are discussed in \cref{sec:surrogate_construction}.  Convergence and performance experiments can be found in \cref{sec:experiments}.  Finally, \cref{sec:conclusion} summarizes the results of the paper and proposes future extensions.


\section{Multirate infinitesimal general-structure additive Runge-Kutta me\-thods}
\label{sec:mri_gark}

In this section, we will briefly review the MRI-GARK framework, as it serves as the foundation for the remainder of the paper.  These methods numerically solve ODEs that are additively partitioned into fast and slow dynamics:
\begin{equation} \label{eq:mr_ode}
	\y' = \fun(t, \y) = \fun\F(t, \y) + \fun\S(t, \y).
\end{equation}

The defining characteristic of multirate infinitesimal methods is that the slow dynamics $\fun\S$ is propagated with a discrete method, most commonly a Runge--Kutta method, while the fast dynamics $\fun\F$ is propagated continuously through modified fast ODEs.  An MRI-GARK scheme \cite{Sandu_2019_MRI-GARK} advances the solution of \eqref{eq:mr_ode} from time $t_n$ to $t_{n+1} = t_n + H$ via the following computational process:
\begin{subequations}
	\label{eq:MRI-GARK-decoupled}
	\begin{align}
		\label{eq:MRI-GARK-decoupled-first_stage}
		& Y\S_1 = \y_n, \\
		\label{eq:MRI-GARK-decoupled-internal_ode}
		& \left\{ \begin{aligned}
			v_i(0) &= Y\S_{i}, \\
			T_{i} &= t_n + \cbase_{i} H, \\
			v_i'(\theta) &= \deltac_{i} \fun\F \bigl(T_{i} + \deltac_{i} \theta, v_i(\theta) \bigr) + \sum_{j=1}^{i+1} \gamma_{i,j} \mleft( \tfrac{\theta}{H} \mright) \fun\S \bigl( T_j,Y\S_j \bigr), \\
			& \quad \text{for } \theta \in [0, H], \\
			Y\S_{i+1} &= v_i(H), \qquad i = 1, \dots, s\S,
		\end{aligned} \right. \\
		\label{eq:MRI-GARK-decoupled-solution}
		& \y_{n+1} = Y\S_{s\S + 1}.
	\end{align}
\end{subequations}
There are $s\S$ modified fast ODEs \cref{eq:MRI-GARK-decoupled-internal_ode} solved between the abscissae of a Runge--Kutta method.  The distances between consecutive abscissae are
\begin{equation*}
	\deltac_{i} = \begin{cases}
		\cbase_{i+1} - \cbase_i, & i = 1, \dots, s\S - 1, \\
		1 - \cbase_{s\S}, & i = s\S.
	\end{cases}
\end{equation*}
\Cref{eq:MRI-GARK-decoupled-internal_ode} also contains time-dependent, polynomial forcing terms for the slow dynamics which satisfy
\begin{equation} \label{eq:gamma_polys}
	\gamma_{i,j}(t) \coloneqq \sum_{k \ge 0} \gamma_{i,j}^k t^k,
	\quad
	\widetilde{\gamma}_{i,j}(t) \coloneqq \int_{0}^{t} \gamma_{i,j}(\tau) \dd \tau = \sum_{k \ge 0} \gamma_{i,j}^k \frac{t^{k+1}}{k+1},
	\quad
	\overline{\gamma}_{i,j} \coloneqq \widetilde{\gamma}_{i,j}(1).
\end{equation}
Continuing with the notation of \cite{Sandu_2019_MRI-GARK}, capitalized versions of \cref{eq:gamma_polys} denote matrices of coefficients.  If $\overline{\Gamma}$ is lower triangular, \cref{eq:MRI-GARK-decoupled} is explicit, and otherwise, it is implicit.

For general $\fun\F$, one cannot expect to find a closed-form solution to the modified fast ODE \cref{eq:MRI-GARK-decoupled-internal_ode}; however, one can determine $v_i(H)$ by time-stepping with ``infinitesimally'' small timesteps.  For the sake of simplifying the analysis, we assume this integration is exact, but computationally, it is not possible to take infinitely many steps.  In practice, \cref{eq:MRI-GARK-decoupled-internal_ode} is integrated to an accuracy that is negligible compared to the other source of error: treating $f\S$ with a discrete Runge--Kutta method.  This is analogous to assuming integrations within a Strang splitting are exact, but using sufficiently accurate approximations in practice.  Multirate infinitesimal methods enjoy great flexibility in that \textit{any} consistent time integration method can be used to solve the modified fast ODEs.  This offers the freedom to choose implicit or explicit methods, as well as methods of any order.  \Cref{eq:MRI-GARK-decoupled} is multirate in the sense that there are $s\S$ evaluations of $\fun\S$ over a timestep $H$, while \cref{eq:MRI-GARK-decoupled-internal_ode} is generally integrated with a smaller timestep and requires many more evaluations of $\fun\F$.

\subsection{Coupled MRI-GARK schemes}

In \cite{Sandu_2020_MRI-GARK_Coupled} two new types of MRI-GARK methods were proposed: step predictor-corrector MRI-GARK (SPC-MRI-GARK) and internal predictor-corrector MRI-GARK (IPC-MRI-GARK).  In order to improve the linear stability of MRI-GARK, both include coupled, discrete prediction stages that are subsequently corrected by solving modified fast ODEs.  In comparing the two MRI-GARK variants, SPC-MRI-GARK was identified to have simpler order conditions, better stability, and a simpler implementation.  For these reasons, we will only consider SPC-MRI-GARK in the construction of surrogate model timestepping techniques.

An SPC-MRI-GARK method based on the ``slow'' Runge--Kutta method $(\Abase, \allowbreak\bbase,\allowbreak \cbase)$ solves \cref{eq:mr_ode} with steps of the form:
\begin{subequations}
	\label{eq:SPC-MRI-GARK}
	\begin{align}
		\label{eq:SPC-MRI-GARK_predictor}
		& Y_i = \y_n + H \sum_{j=1}^{s\S} \abase_{i,j} \fun \bigl( t_n + \cbase_{j} H, Y_j \bigr), \qquad i =1,\dots,s\S, \\
		\label{eq:SPC-MRI-GARK_corrector}
		& \left\{ \begin{aligned}
			v(0) &= \y_n, \\
			v' &= \fun\F(t_n + \theta, v) + \sum_{j=1}^{s\S} \gamma_j \mleft( \tfrac{\theta}{H} \mright) \fun\S \bigl( t_n + \cbase_{j} H, Y_j \bigr), \quad \text{for } \theta \in [0, H], \\
			y_{n+1} &= v(H).
		\end{aligned} \right.
	\end{align}
\end{subequations}

Similar to a traditional Runge--Kutta scheme, SPC-MRI-GARK starts by computing stage values that approximate the solution at the times $t_n + \cbase_{i} H$.  In this prediction step \cref{eq:SPC-MRI-GARK_predictor}, the fast and slow dynamics are treated together with the timestep $H$.  The fast part of the stage values is inaccurate from this large timestep and is discarded.  A single ODE \cref{eq:SPC-MRI-GARK_corrector} is solved from $t_n$ to $t_{n+1}$ to correct the fast dynamics and produce the next step $y_{n+1}$.  Like MRI-GARK, there are time-dependent, polynomial forcing terms for the slow tendencies.  Following \cite[Definition 2.2]{Sandu_2020_MRI-GARK_Coupled}, the definitions of $\gamma_{i}(t)$, $\widetilde{\gamma}_{i}(t)$, $\overline{\gamma}_{i}$ are identical to \cref{eq:gamma_polys}, except they are vector quantities with a single subscript index.
While base methods of any structure are admitted in the SPC-MRI-GARK framework, only diagonally implicit methods were considered in \cite{Sandu_2020_MRI-GARK_Coupled}.  As the focus of this paper is explicit methods, we have derived new explicit methods of orders two and three and present their coefficients in \cref{app:method_coeffs}.

\section{Method formulation}
\label{sec:method_formulation}

Following the process outlined in the introduction, we seek to apply the MRI-GARK method \cref{eq:MRI-GARK-decoupled} to the partitioned problem \cref{eq:full_surrogate_ode} in order to construct a solution process that can incorporate surrogate model information.  We start with the modified fast ODE \cref{eq:MRI-GARK-decoupled-internal_ode}:
\begin{equation} 
\label{eq:modified_ode_no_simplify}
\begin{split}
	v_i'(\theta) &= \deltac_{i} \Vmat \fsur \bigl(T_{i} + \deltac_{i} \theta, \Wmat* v_i(\theta) \bigr) \\
	&\quad + \sum_{j=1}^{i+1} \gamma_{i,j} \mleft( \tfrac{\theta}{H} \mright) \bigl( \fun(T_j, Y_j) - \Vmat \fsur \mleft(T_j, \Wmat* Y_j \mright) \bigr).
	\end{split}
\end{equation}
In \eqref{eq:modified_ode_no_simplify}, the full model $\fun$ only appears in the forcing terms, and for explicit methods, it is evaluated at previously computed stage values.  Only the surrogate model $\fsur$ needs to be evaluated at each of the infinitesimal steps to compute $v_i(H)$.

Note that in formulation \eqref{eq:modified_ode_no_simplify}, $v_i \in \Cplx{\nvar}$ evolves in the full space, and there is no benefit from choosing a surrogate model that evolves in a lower dimensional space.  To resolve this, we first split $v_i$ into the parts residing inside and outside the surrogate model space:
\begin{equation*}
	v_i(\theta) = \Vmat \Wmat* v_i(\theta) + \left( \eye{\nvar} - \Vmat \Wmat* \right) v_i(\theta).
\end{equation*}
More formally, $\Vmat \Wmat*$ and $\left( \eye{\nvar} - \Vmat \Wmat* \right)$ define oblique projections onto the range of $\Vmat$ and the nullspace of $\Wmat*$, respectively.  We now consider the solution to \cref{eq:modified_ode_no_simplify} from the perspective of these two, complementary subspaces.

Starting outside the surrogate model space, we have that
\begin{equation} \label{eq:ode_solution_outside}
	\begin{split}
		\left( \eye{\nvar} - \Vmat \Wmat* \right) v_i(H)
		&= \left( \eye{\nvar} - \Vmat \Wmat* \right) \left( Y_i + \int_{0}^{H} v_i'(\theta) \dd{\theta} \right) \\
		&= \left( \eye{\nvar} - \Vmat \Wmat* \right) \left( Y_i + H \sum_{j=1}^{i} \overline{\gamma}_{i,j} \fun(T_j, Y_j) \right),
	\end{split}
\end{equation}
where we have used \cref{eq:gamma_polys} and the fact that $\left( \eye{\nvar} - \Vmat \Wmat* \right) \Vmat = \zero{\nvar}{\nsur}$ by \cref{eq:space_transformation}.  \Cref{eq:ode_solution_outside} reveals that inside the nullspace of $\Wmat*$ the solution to equation \cref{eq:modified_ode_no_simplify} becomes a Runge--Kutta stage.  Therefore, we exclude the components in this subspace during the infinitesimal step integration of \cref{eq:modified_ode_no_simplify}.

\begin{figure}
	\centering
	\includegraphics[width=0.65\linewidth]{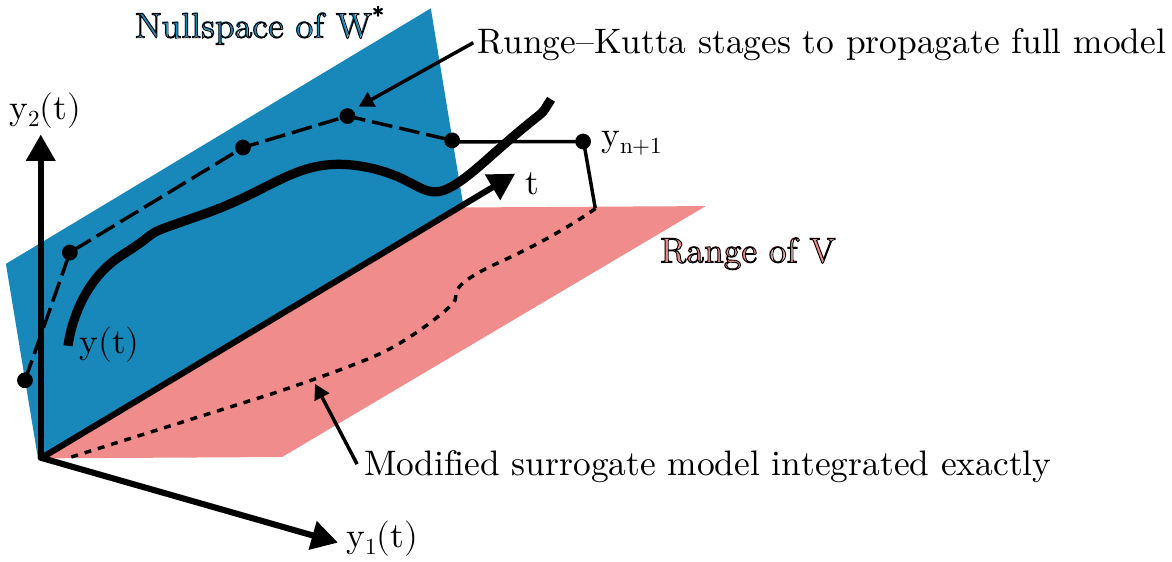}
	\caption{Illustration of \mrimethod{} for a two-variable ODE and a surrogate model that evolves in a one-dimensional subspace.}
	\label{fig:SM-MRI-GARK}
\end{figure}

We now consider the dynamics \cref{eq:modified_ode_no_simplify} projected onto the surrogate model space:
\begin{equation} \label{eq:ode_solution_inside}
	\begin{split}
		z_i(\theta) &\coloneqq \Wmat* v_i(\theta)  \in \Cplx{\nsur}, \\
		z_i'(\theta) &= \deltac_i \fsur \bigl( T_i + \deltac_i \theta, z_i(\theta) \bigr) \\
		& \quad + \sum_{j=1}^{i+1} \gamma_{i,j} \mleft( \tfrac{\theta}{H} \mright) \left( \Wmat* \fun(T_j, Y_j) - \fsur \mleft(T_j, \Wmat* Y_j \mright) \right).
	\end{split}
\end{equation}
The solution to \cref{eq:modified_ode_no_simplify} can be expressed as the sum of its two components:
\begin{equation*}
	v_i(H) = \Vmat z_i(H) + \left( \eye{\nvar} - \Vmat \Wmat* \right) \left( Y_i + H \sum_{j=1}^{i+1} \overline{\gamma}_{i,j} \fun(T_j, Y_j) \right).
\end{equation*}
Thus to compute an MRI-GARK stage, we only need to solve an ODE of dimension $\nsur$ instead of $\nvar$.  We summarize the simplifications of \cref{eq:MRI-GARK-decoupled} in \cref{def:SM-MRI-GARK} and provide a graphical illustration in \cref{fig:SM-MRI-GARK}.

\begin{definition}[\mrimethod{}] \label{def:SM-MRI-GARK}
	A surrogate model MRI-GARK (\mrimethod{}) method solves the ODE \cref{eq:ode} with the help of the surrogate model \cref{eq:surrogate_ode} using steps of the form
	\begin{subequations} \label{eq:SM-MRI-GARK}
		\begin{align}
			& Y_1 = \y_n, \label{eq:SM-MRI-GARK:Y1} \\
			& \left\{
			\begin{aligned}
				z_i(0) &= \Wmat* Y_i, \\
				T_i &= t_n + \cbase_i H, \\
				z_i'(\theta) &= \deltac_i \fsur \bigl( T_i + \deltac_i \theta, z_i(\theta) \bigr) \\
				& \quad + \sum_{j=1}^{i+1} \gamma_{i,j} \mleft( \tfrac{\theta}{H} \mright) \left( \Wmat* \fun(T_j, Y_j) - \fsur \mleft(T_j, \Wmat* Y_j \mright) \right), \quad \text{for } \theta \in [0, H], \\
				Y_{i+1} &= \Vmat z_i(H) + \left( \eye{d} - \Vmat \Wmat* \right) \left( Y_i + H \sum_{j=1}^{i+1} \overline{\gamma}_{i,j} \fun(T_j, Y_j) \right), \\
				& \quad i = 1, \dots, s\S,
			\end{aligned}
			\right. \label{eq:SM-MRI-GARK:ode} \\
			& \y_{n+1} = Y_{s\S + 1}. \label{eq:SM-MRI-GARK:yn}
		\end{align}
	\end{subequations}
\end{definition}

\begin{remark}[\mrimethod{} for hierarchical surrogate models]
	For certain applications, an entire hierarchy of surrogate models $\fun, \fsur^{[1]}, \fsur^{[2]}, \dots, \fsur^{[m]}$ may be available.  In such scenarios, one can take advantage of all models by applying \mrimethod{} recursively.  First, the highest fidelity models $\fun$ and $\fsur^{[1]}$ are applied to \cref{eq:SM-MRI-GARK}.  The ODE \cref{eq:SM-MRI-GARK:ode}, is then solved using an \mrimethod{} method with $\deltac_i \fsur^{[2]}$ as the surrogate model.  Its ODEs are solved using an \mrimethod{} method with $\deltac[\times 2]_i \fsur^{[3]}$ as the surrogate model.  This is continued until we reach $\deltac[\times (m-1)]_i \fsur^{[m]}$ as the lowest fidelity surrogate model.  Note the power of $\deltac_i$ scaling the surrogate models can be eliminated by applying the change of variable $\theta \to \theta / \deltac_i$ to \cref{eq:SM-MRI-GARK:ode}.
\end{remark}

\begin{remark}[Connection to EPIRK-K methods]
	In \cite{Sandu_2019_MRI-GARK}, Sandu showed that exponential Runge--Kutta methods are closely related to MRI-GARK methods when the linear-nonlinear partitioning 
	\begin{equation*}
		\y' = \underbrace{\Lb\,\y}_{\fun\F(t, \y)} + \underbrace{\fun(t, \y) - \Lb\,\y}_{\fun\S(t, \y)}
	\end{equation*}
	is used.  MRI-GARK can be viewed as a generalization of exponential Runge--Kutta methods where the piece integrated exactly can be linear or nonlinear.
	
	An analogous connection can be made between EPIRK-K \cite{Sandu_2019_EPIRKW} and \mrimethod{}.  EPIRK-K is an extension of exponential propagation iterative methods of Runge-Kutta type (EPIRK) \cite{tokman2011new} that approximates the Jacobian $\Jac_n = {\pdv{\fun}{\y}}(t_n, \y_n)$ appearing in matrix exponentials and related matrix functions with a Krylov-subspace approximation $\mathbf{U} \mathbf{U}^* \Jac_n \mathbf{U} \mathbf{U}^*$.  \mrimethod{} with $\Vmat = \Wmat = \mathbf{U}$ and $\fsur(t, \ysur) = \mathbf{U}^* \Jac_n \mathbf{U} \ysur$ closely resembles the class of EPIRK-K schemes.
\end{remark}

\begin{remark}[Stepsize adaptivity]
	\label{rem:adaptivity}
	Even with fixed timesteps $H$, adaptive time-steppers can be used to integrate \cref{eq:SM-MRI-GARK:ode} so that $\fsur$ is evaluated just enough to meet tolerances.  It is also possible to independently and adaptively select $H$, and thus, the frequency at which $\fun$ is evaluated.  
	Following \cite[Remark 2.5]{Sandu_2019_MRI-GARK}, \cref{eq:SM-MRI-GARK} can be equipped with an embedded method to estimate the local truncation error, and then standard error controllers can be used \cite[Section II.4]{Hairer1993}.  Each method in \cref{app:method_coeffs} includes an embedded method that is one order lower than the main method.
\end{remark}

As a simple example of an \mrimethod{} method, consider forward Euler as the base method and the internally consistent \cite[Section 3.1.1]{Sandu_2019_MRI-GARK} coupling $\gamma_{1,1}(t) = 1$:
\begin{equation} \label{eq:SM-MRI-GARK_Euler}
	\begin{split}
		z(0) &= \Wmat* \y_n \\
		z'(\theta) &= \fsur(t_n + \theta, z(\theta)) + \Wmat* \fun(t_n, \y_n) - \fsur \mleft(t_n, \Wmat* \y_n \mright), \qquad \text{for } \theta \in [0, H], \\
		\y_{n+1} &= \Vmat z(H) + \left( \eye{d} - \Vmat \Wmat* \right) \left( \y_n + H \fun(t_n, \y_n) \right).
	\end{split}
\end{equation}
\Cref{eq:SM-MRI-GARK_Euler} is only first order accurate; however, second and third order \mrimethod{} schemes can be constructed using the coefficents provided in \cref{app:method_coeffs}.  In order to implement these efficiently, we have provided pseudocode in \cref{alg:SM-MRI-GARK}.  It minimizes the number of full and surrogate model evaluations, as well as the number of matrix-vector products involving $\Vmat$ and $\Wmat*$.

\begin{algorithm}[ht!]
	\caption{Pseudocode for explicit \mrimethod{} \cref{eq:SM-MRI-GARK}.}
	\label{alg:SM-MRI-GARK}
	\begin{algorithmic}[1]
		\Procedure{\mrimethod}{$\fun, \fsur, \Vmat, \Wmat*, t_0, t_\text{end}, \y_0, N_\text{steps}$}
			\State $y = \y_0$
			\State $\widehat{y} = \Wmat* \y_0$
			\State $H = (t_\text{end} - t_0) / N_\text{steps}$
			
			\For{$n = 1, \dots, N_\text{steps}$}
				\For{$i = 1, \dots, s\S$}
					\State $T_i = t_0 + ( n + \cbase_i ) H$
					\State $k_i = \fun(T_i, y)$ 
					\State $\widehat{k}_i = \Wmat* k_i$
					\State $\widehat{\ell}_i = \widehat{k}_i - \fsur(T_i, \widehat{y})$
					\State $\widehat{w} = \widehat{y} + H \sum_{j=1}^{i} \overline{\gamma}_{i,j} \widehat{k}_j$
					\State $y = y + H \sum_{j=1}^{i} \overline{\gamma}_{i,j} k_j$
					\State $\widehat{y} = \text{OdeSolve}\mleft( z_i'(\theta) = \deltac_i \fsur ( T_i + \deltac_i \theta, z_i(\theta) ) + \sum\limits_{j=1}^{i} \gamma_{i,j} \mleft( \tfrac{\theta}{H} \mright) \widehat{\ell}_j, \mright.$ \newline
					\hspace*{10.75em} $\mleft. \text{timespan} = [0, H], \text{initial\_condition} = \widehat{y} \mright)$
					\State $y = y + \Vmat (\widehat{y} - \widehat{w})$
				\EndFor
			\EndFor
			
			\State \textbf{return} $y$
		\EndProcedure
	\end{algorithmic}
\end{algorithm}

\subsection{Formulation for SPC-MRI-GARK}

A similar process can be used to apply SPC-MRI-GARK \cref{eq:SPC-MRI-GARK} to the ODE \cref{eq:full_surrogate_ode}.  For brevity, we avoid replicating the simplifications \cref{eq:ode_solution_outside,eq:ode_solution_inside} and directly present the method formulation in the following definition.

\begin{definition}[\spcmethod{}] \label{def:SM-SPC-MRI-GARK}
	A surrogate model SPC-MRI-GARK (\spcmethod{}) method solves the ODE \cref{eq:ode} with the help of the surrogate model \cref{eq:surrogate_ode} using steps of the form
	\begin{subequations} \label{eq:SM-SPC-MRI-GARK}
		\begin{align}
			& Y_i = \y_n + H \sum_{j=1}^{s\S} \abase_{i,j} \fun \bigl( t_n + \cbase_{j} H, Y_j \bigr), \qquad i =1,\dots,s\S, \label{eq:SM-SPC-MRI-GARK:stage} \\
			& \left\{
			\begin{aligned}
				z(0) &= \Wmat* \y_n,\\
				z(\theta)' &= \fsur \mleft( t_n + \theta, z(\theta) \mright) \\
				& \quad + \sum_{j=1}^{s\S} \gamma_{j} \mleft( \tfrac{\theta}{H} \mright) \left( \Wmat* \fun \bigl( t_n + \cbase_{j} H, Y_j \bigr) - \fsur \bigl( t_n + \cbase_{j} H, \Wmat* Y_j \bigr) \right), \\
				& \quad \text{for } \theta \in [0, H],\\
				y_{n+1} &= \Vmat z(H) + \left( \eye{d} - \Vmat \Wmat* \right) \left( \y_n + H \sum_{j=1}^{s\S} b_{j} \fun \bigl( t_n + \cbase_{j} H, Y_j \bigr) \right).
			\end{aligned}
			\right. \label{eq:SM-SPC-MRI-GARK:ode}
		\end{align}
	\end{subequations}
\end{definition}

The \mrimethod{} Euler method from \cref{eq:SM-MRI-GARK_Euler} also happens to be an \spcmethod{} method.  The base method is again forward Euler, and the coupling is $\gamma_{1}(t) = 1$.  Methods with more than one stage will not coincide, however, as \spcmethod{} has one modified fast ODE per step while \mrimethod{} has $s\S$.  Another distinction, which can be seen in \cref{alg:SM-SPC-MRI-GARK}, is that \spcmethod{} only requires one matrix-vector product with $\Vmat$ per step, while \mrimethod{} requires $s\S$.

\begin{algorithm}[ht!]
	\caption{Pseudocode for explicit \spcmethod{} \cref{eq:SM-SPC-MRI-GARK}.}
	\label{alg:SM-SPC-MRI-GARK}
	\begin{algorithmic}[1]
		\Procedure{\spcmethod}{$\fun, \fsur, \Vmat, \Wmat*, t_0, t_\text{end}, \y_0, N_\text{steps}$}
		\State $y = \y_0$
		\State $\widehat{y} = \Wmat* \y_0$
		\State $H = (t_\text{end} - t_0) / N_\text{steps}$
		
		\For{$n = 1, \dots, N_\text{steps}$}
			\State $t_n = t_0 + n H$
			\For{$i = 1, \dots, s\S$}
				\State $k_i = f \bigl( t_n + \cbase_i H, y + H \sum_{j=1}^{i-1} \abase_{i,j} k_j \bigr)$ 
				\State $\widehat{k}_i = \Wmat* k_i$
				\State $\widehat{\ell}_i = \widehat{k}_i - \fsur \bigl( t_n + \cbase_i H, \widehat{y} + H \sum_{j=1}^{i-1} \abase_{i,j} \widehat{k}_j \bigr)$
			\EndFor
			
			\State $\widehat{w} = \widehat{y} + H \sum_{j=1}^{s\S} \bbase_{j} \widehat{k}_j$
			\State $y = y + H \sum_{j=1}^{s\S} \bbase_{j} k_j$
			\State $\widehat{y} = \text{OdeSolve}\mleft( z'(\theta) = \fsur \mleft( t_n + \theta, z(\theta) \mright) + \sum_{j=1}^{s\S} \gamma_{j} \mleft( \tfrac{\theta}{H} \mright) \widehat{\ell}_j, \mright.$ \newline
			\hspace*{9.25em} $\mleft. \text{timespan} = [0, H], \text{initial\_condition} = \widehat{y} \mright)$
			\State $y = y + \Vmat (\widehat{y} - \widehat{w})$
		\EndFor
		
		\State \textbf{return} $y$
		\EndProcedure
	\end{algorithmic}
\end{algorithm}


\section{Error analysis}
\label{sec:error}

Consider the edge case where the surrogate model is taken to be the full model and $\Vmat = \Wmat* = \eye{\nvar}$.  The stages of \cref{eq:SM-MRI-GARK} simplify to
\begin{align*}
	& Y_1 = \y_n, \\
	& \left\{
	\begin{aligned}
		z_i(0) &= Y_i, \\
		z_i' &= \deltac_i \fsur \bigl( t_n + \cbase_i H + \deltac_i \theta, z_i \bigr), \qquad \text{for } \theta \in [0, H],\\
		Y_{i+1} &= z_i(H), \qquad i = 1, \dots, s\S,
	\end{aligned}
	\right. \\
	& \y_{n+1} = Y_{s\S + 1}.
\end{align*}
This is simply the original ODE \cref{eq:ode} solved between the abscissae of a Runge--Kutta method.  With the infinitesimal step assumption, the error is zero.  At the other extreme, take the surrogate model to be $\fsur(t, \ysur) = \zero{\nvar}$.  Again, \cref{eq:SM-MRI-GARK} can be simplified, but now into a traditional, unpartitioned Runge--Kutta method:
\begin{align*}
	Y_1 &= \y_n, \\
	Y_{i+1} &= Y_i + H \sum_{j=1}^{i} \overline{\gamma}_{i,j} \fun \bigl( t_n + \cbase_i H, Y_j \bigr), \qquad i = 1, \dots, s\S, \\
	\y_{n+1} &= Y_{s\S + 1}.
\end{align*}
This error will be nonzero in general and can be analyzed using classical Runge--Kutta order condition and convergence theory.

While neither of these two choices for the surrogate model are particularly interesting, they illustrate how the accuracy of \mrimethod{} and \spcmethod{} depends on several factors including the accuracy of the surrogate model, the projection matrices $\Vmat$ and $\Wmat$, and the method coefficients.  In this section, we analyze the structure of local truncation error
\begin{equation}
	\label{eq:lte}
	\textrm{LTE}_{n+1} = \y(t_{n+1}) - \y_{n+1} = C\, H^{p+1} + \order{H^{p+2}}
\end{equation}
to better understand the influence of these factors.  To simplify our error analysis, we use the autonomous form of \cref{eq:full_surrogate_ode}, which comes at no loss of generality for internally consistent methods.

The new families of methods in this paper are derived as special cases of MRI-GARK and SPC-MRI-GARK.  Interestingly, it is also possible to view MRI-GARK and SPC-MRI-GARK as special cases of the surrogate model methods.  When
\begin{equation*}
	\Vmat = \Wmat* = \eye{\nvar},
	\quad
	\fun(t, \y) = \fun\F(t, \y) + \fun\S(t, \y),
	\quad
	\fsur(t, \y) = \fun\F(t, \y),
\end{equation*}
we identically recover \cref{eq:MRI-GARK-decoupled,eq:SPC-MRI-GARK}.  Thanks to these properties, MRI-GARK and SPC-MRI-GARK order conditions are both necessary and sufficient for the surrogate model versions.  Interested readers can refer to \cite[Section 3.1]{Sandu_2019_MRI-GARK} and \cite[Section 2.2]{Sandu_2020_MRI-GARK_Coupled} for these order conditions, which are derived using N-tree theory \cite{Araujo1997}.  Thus, the order of the surrogate-model-based integration only depends on the coefficients of the underlying MRI-GARK method.  Notably, the order does not depend on the accuracy or quality of the surrogate model. 

We note that these results are based on nonstiff order condition theory which requires $\fun$ and $\fsur$ to be $p+1$ times continuously differentiable and have moderate Lipshitz constants with respect to $H$.  In cases where these assumptions do not hold, e.g., nonsmooth surrogate models or stiff, infinite-dimensional PDEs, order reduction can occur.  

We now discuss the leading error constant $C$ in \cref{eq:lte}.  N-Trees provide an elegant way to quantify the terms of which it is composed.  For our multirate methods, there are two partitions, so we use the set of two-trees $\BT_2$.  Tree vertices are ``colored'' either fast (black circle) or slow (white circle).  Let $\BT_2\F$ be the set of trees in $\BT_2$ with all vertices fast-colored (including the empty tree), and $\BT_2 \setminus \BT\F_2$  the set of trees with at least one slow-colored node. We have that
\begin{equation} \label{eq:error_const}
	\begin{split}
		C &= \sum_{\substack{\utree \in \BT_2 \setminus \BT\F_2 \\ \rho(\utree) = p+1}} \frac{1}{\sigma(\utree)} \left( \Phi(\utree) - \frac{1}{\gamma(\utree)} \right) F(\utree)(y), 
		\\
		 \BT_2 \setminus \BT\F_2 &= 
		\Bigl\{
		\raisebox{-5pt}{$
			\btree{[s]},
			\btree{[s[s]]}, \btree{[s[f]]}, \btree{[f[s]]},
			\btree{[s[s][s]]}, \btree{[s[f][s]]}, \btree{[s[s][f]]}, \btree{[s[f][f]]}, \btree{[f[s][s]]}, \btree{[f[f][s]]}, \btree{[f[s][f]]},
			\btree{[s[s[s]]]},
			\dots
		$}
		\Bigr\},
	\end{split}
\end{equation}
where $\rho(\utree)$, $\sigma(\utree)$, $\gamma(\utree)$, are the order, symmetries, and density, respectively \cite{Araujo1997}.  The elementary weight $\Phi(\utree)$ depends only on the method and its coefficients.  For $\utree \in \BT\F_2$, $\Phi(\utree) = 1 / \gamma(\utree)$ by the infinitesimal step assumption, which allows us to exclude this subset of trees from the summation in \cref{eq:error_const}.  Finally, the elementary differentials have the recursive, tensorial definition:
\begin{subequations}
	\begin{align}
		F(\emptyset)(\y) &= \y, \\
		F \bigl( \btree{[f]} \bigr)(\y) &= \Vmat \fsur(\Wmat* \y), \\
		F \bigl( \btree{[s]} \bigr)(\y) &= \fun(\y) - \Vmat \fsur(\Wmat* \y), \\
		F \mleft( \left[ \utree_1, \dots, \utree_m \right]\F \mright)(\y) &= \Vmat \fsur^{(m)}( \Wmat* \y ) \bigl( \Wmat* F(\utree_1)(\y), \dots, \Wmat* F(\utree_m)(\y) \bigr), \\
		F \mleft( \left[ \utree_1, \dots, \utree_m \right]\S \mright)(\y) &= \fun^{(m)}(\y) \bigl( F(\utree_1)(\y), \dots, F(\utree_m)(\y) \bigr) - F \mleft( \left[ \utree_1, \dots, \utree_m \right]\F \mright)(\y). \label{eq:elem_diff_s}
	\end{align}
\end{subequations}
A tree $\utree \in \BT_2$ can be expressed as $\left[ \utree_1, \dots, \utree_m \right]\comp{\sigma}$ where $\sigma$ is the color of the root vertex and $\utree_1, \dots, \utree_m$ are the nonempty subtrees left when the root is removed.  \Cref{tab:trees} provides examples of terms appearing in \cref{eq:error_const}.

\begin{table}[ht!]
	\centering
	\begin{tabular}{c|cccc}
		$\utree$ & $F(\utree)(y)$ & $\rho(\utree)$ & $\sigma(\utree)$ & $\gamma(\utree)$ \\ \hline
		$\btree{[s]}$ & $\fun(\y) - \Vmat \fsur(\Wmat* \y)$ & 1 & 1 & 1 \\
		$\btree{[s[f]]}$ & $\fun'(\y) \Vmat \fsur(\Wmat* \y) - \Vmat \fsur'(\Wmat* \y) \fsur(\Wmat* y)$ & 2 & 1 & 2 \\
		$\btree{[f[s][f]]}$ & $\Vmat \fsur''(\Wmat* \y) \bigl( \Vmat \fsur(\Wmat* \y), \fun(\y) - \Vmat \fsur(\Wmat* \y) \bigr)$ & 3 & 1 & 3
	\end{tabular}
	\caption{Examples of trees, elementary differentials, and other tree properties.}
	\label{tab:trees}
\end{table}

Using \cref{eq:error_const}, the local truncation error of \mrimethod{} Euler from \cref{eq:SM-MRI-GARK_Euler} is
%
\begin{equation*}
	\textrm{LTE}_{n+1}^{\text{Euler}} = \frac{1}{2} ( \fun'(\y_n) - \Vmat \fsur'(\Wmat* \y_n) \Wmat* ) \fun(\y_n) H^2 + \order{H^3}.
\end{equation*}

Note that each tree in the summation in \cref{eq:error_const} contains at least one slow node.  From \cref{eq:elem_diff_s}, a slow-colored vertex corresponds to the surrogate model error or one of its derivatives appearing in an elementary differential.  In order to quantify the error in the surrogate model, suppose there exists a (small) constant $\varepsilon$ such that:
\begin{align*}
	\norm{\fun(\y_n) - \Vmat \fsur(\Wmat* \y_n)} &\leq \varepsilon, \\
	\norm{\fun^{(m)}(\y_n) \mleft( x_1, \dots, x_m \mright) - \Vmat \fsur^{(m)}( \Wmat* \y_n ) \mleft( \Wmat* x_1, \dots, \Wmat* x_m \mright)}
	&\leq \varepsilon \norm{x_1} \cdots \norm{x_m},
\end{align*}
for $m = 1, \dots p$.  Here $\norm{\cdot}$ denotes an arbitrary norm on $\Cplx{\nvar}$.  This yields the bounds
\begin{alignat*}{3}
	\norm{F \bigl( \btree{[f]} \bigr)(\y_n)} &\leq d_0 + \varepsilon,
	\qquad &
	\norm{F \mleft( \left[ \utree_1, \dots, \utree_m \right]\F \mright)(\y_n)} &\leq (d_m + \varepsilon) \prod_{i=1}^m \norm{F(u_i)(\y_n)}, \\
	\norm{F \bigl( \btree{[s]} \bigr)(\y_n)} &\leq \varepsilon,
	\qquad &
	\norm{F \mleft( \left[ \utree_1, \dots, \utree_m \right]\S \mright)(\y_n)} &\leq \varepsilon \prod_{i=1}^m \norm{F(u_i)(\y_n)},
\end{alignat*}
where $d_0 = \norm{\fun(\y_n)}$ and $d_m$ is the operator norm of $\fun^{(m)}(\y_n)$.  Now we have that
\begin{equation*}
	\norm{C} \leq \sum_{\substack{\utree \in \BT_2 \setminus \BT\F_2 \\ \rho(\utree) = p+1}} \abs{\frac{1}{\sigma(\utree)} \left( \Phi(\utree) - \frac{1}{\gamma(\utree)} \right)} \varepsilon^{\rho\S(\utree)} \left( \varepsilon + \max_{i = 0, \dots, p} d_i \right)^{\rho\F(\utree)} = \sum_{i=1}^{p+1} c_i \varepsilon^i,
\end{equation*}
where $\rho\comp{\sigma}(\utree)$ is the number of $\sigma$-colored vertices in $\utree$.  The constants $c_i$ only depend on $\fun$ and the order condition residuals; they are independent of $\fsur$, $\Vmat$, $\Wmat$, $H$ and $\varepsilon$.  As one might expect, the local truncation error decreases as the accuracy of the surrogate model and its derivatives increase.

%
%

\section{Construction of surrogate models for accelerating time integration}
\label{sec:surrogate_construction}

This sections discusses several techniques to construct a surrogate model $\fsur$ and the associated linear operators $\Vmat$ and $\Wmat*$ with a computationally favorable balance of accuracy and evaluation cost.

\subsection{Reduced-order models (ROMs)}

There are numerous techniques from the reduced-order modeling community that may be used to generate the surrogate models including proper orthogonal decomposition (POD) \cite{sirovich1987turbulence1}, Krylov-subspace methods \cite{grimme1997krylov}, reduced basis methods \cite{peterson1989reduced}, discrete empirical interpolation method \cite{chaturantabut2010nonlinear}, and dynamic mode decomposition \cite{schmid2010dynamic,tissot2014model}.
Other multi-resolution methods based on Fourier or wavelet transformation \cite{brunton2019data} could be used to only capture coarse information about the model.

\subsection{Multimesh models}
\label{subsec:multimesh}

Consider a numerical approach to discretize a partial differential equation (PDE) in space over a computational mesh using, for example, the finite element, difference, or volume method.  A surrogate model can come from solving the PDE on a coarser mesh or using a lower order spatial discretization.  This provides an approximation to capture the ``shape'' of the solution that is cheaper and also enjoys a less strict Courant--Friedrichs--Lewy (CFL) condition.  If the meshes are nested, $\Wmat$ is a subset of columns of the identity matrix, and $\Vmat$ is a sparse interpolation matrix.  We note that this strategy closely resembles multigrid methods\cite{zubair2009efficient}.  The relaxation and prolongation operators, however, have to be chosen such as  $\Wmat* \Vmat=\eye{\nvar}$,
such as the 1D relaxation stencil $\begin{bsmallmatrix} 0 & 1 & 0\end{bsmallmatrix}$ and the 1D prolongation stencil $\begin{bsmallmatrix} 1/2 & 1 & 1/2\end{bsmallmatrix}$.

\subsection{Applying simplifying approximations to the full model}

By linearizing, filtering, averaging, simply ignoring, or otherwise approximating certain terms in $\fun$, a computationally inexpensive $\fsur$ can be produced.  In \cite[Section 6]{chacon2017multiscale}, for example, a surrogate model of z-level ocean model is produced by a vertical averaging of barotropic velocities.  Consider also a direct N-body simulation with a Barnes-Hut or fast multipole method used as the surrogate model.  For this example, the surrogate models happen to have the same state representation as the full model, and therefore, $\nsur = \nvar$ and $\Vmat = \Wmat = \eye{\nvar}$.  


\subsection{Data-driven surrogate models}

When sufficient input and output data is available, supervised learning approaches are a viable option to generate approximate models for a system.  A variety of techniques have been developed, some depending on system identification and sparse dictionary learning \cite{Kutz2017ModelDiscovery}, others using neural networks to discover operators and right-hand side functions \cite{Raissi2019Feb}. Patch data, both in time and space \cite{Kevrekidis2004EquationFree} have been used to train machine learning models that can reproduce crude or partial dynamics of the full system.
In some cases, the data driven dynamics reside in the full space, so $\nsur=\nvar$ and $\Vmat = \Wmat = \eye{\nvar}$.  In other cases, the dynamics could reside in the dominant modes of the data, in which case $\nsur < \nvar$ and $\Vmat$ and $\Wmat$ are determined by the method.


\section{Numerical experiments}
\label{sec:experiments}

In order to evaluate the new surrogate model time-steppers, we apply them to a diverse set of ODE test problems equipped with surrogate models.  Two of the test problems are used to verify convergence properties (error versus steps) and another two are used to measure the integrators' performance (error versus runtime).  The methods considered in the experiments are \mrimethod{} Euler from \cref{eq:SM-MRI-GARK_Euler} and the four methods from \cref{app:method_coeffs}.  We compare these with their base methods, which are traditional Runge--Kutta methods, when only the full or only the surrogate model is used.  In all the experiments, the error is computed by comparing a numerical solution to a high-accuracy reference solution at the final time.  The solution using only the surrogate model resides in the surrogate model subspace, so it is multiplied by $\Vmat$ before computing error.

\subsection{Quasi-geostrophic model and quadratic ROM}
The quasi-geostrophic (QG) equations~\cite{foster2013finite,foster2016conforming,ferguson2008numerical}, are a well-studied set of approximations to both atmospheric and ocean flow. The QG equations have a wide range of well-studied lower dimensional approximations. For our use case, a quadratic POD-Galerkin ROM makes an excellent surrogate model.

We follow the formulations of~\cite{mou2019data,san2015stabilized}, and the same setup as utilized in~\cite{Sandu_2021_MF-EnKF}:
\begin{equation} \label{eq:qge}
	\begin{split}
		\pdv{\omega}{t} + J(\psi,\omega) - \mathrm{Ro}^{-1} \pdv{\psi}{x} &= \mathrm{Re}^{-1}\laplacian{\omega} + \mathrm{Ro}^{-1}F, \\
		\omega &= -\laplacian{\psi}, \\
		J(\psi,\omega) &\equiv \pdv{\psi}{y} \pdv{\omega}{x} - \pdv{\psi}{x} \pdv{\omega}{y}.
	\end{split}
\end{equation}
Here $\omega$ represents the vorticity, $\psi$ the streamfunction, $\mathrm{Re}=450$ the Reynolds number, and $\mathrm{Ro}=0.0036$ the Rossby number. The forcing term is selected to be a symmetric double Gyre to simulate flow in the ocean~\cite{mou2019data,san2015stabilized}:
\begin{equation*}
	F = \sin(\pi(y-1)).
\end{equation*}
The domain for the problem is $\Omega=[0,1]\times[0,2]$. Homogeneous Dirichlet,
\begin{equation*}
	\psi(x,y) = 0,\quad \forall(x,y)\in\partial\Omega,
\end{equation*}
boundary conditions are used for all time.

For the spatial discretization, a second order finite difference is performed, with the canonical Arakawa approximation~\cite{arakawa1966computational} performed on the Jacobian term, $J$. The domain is discretized using 63 points in the $x$ direction and 127 points in the $y$ direction. The discretization begets the discrete stream function variable $\psi$. The relation between the streamfunction and vorticity is used to define the time derivative of the streamfunction variable $\pdv{\psi}{t}$, which we will take as the PDE of interest.

Following~\cite{mou2019data}, and using the method of snapshots~\cite{sirovich1987turbulence1}, we construct basis transformation operators that capture the time-averaged dominant linear modes of the system.
The basis transformations are represented by the orthogonal, in $L^2$, linear operators $\Vmat$, and $\Wmat*$, where $\ysur = \Wmat*\psi$ represents the basis transformation from streamfunction space into ROM space. Details as to how these are derived are found in~\cite{Sandu_2021_MF-EnKF}.

In order to construct a reduced-order model we have to optimally approximate the time derivative of the reduced-order state.
As \cref{eq:qge} is quadratic in nature, we can construct a quadratic POD-Galerkin reduced-order model, of the form
\begin{equation*}
	\ysur' = \fsur(t, \ysur) = \Wmat* \fun(t, \Vmat \ysur) = b + \mathbf{A}\ysur + \ysur^T \mathcal{B} \ysur,
\end{equation*}
where $\ysur$ is the state in the POD basis, $\fun$ is the solution to the Poisson equation in \cref{eq:qge} for the streamfunction derivative $\pdv{\psi}{t}$, $b\in\R{\nsur}$ corresponds to the forcing term $\mathrm{Ro}^{-1}F$, $\mathbf{A}\in\R{\nsur}{\nsur}$ corresponds to the linear term $\mathrm{Ro}^{-1}\pdv{\psi}{x} + \mathrm{Re}^{-1}\laplacian{\omega}$, and $\mathcal{B}\in\R{\nsur}{\nsur}{\nsur}$ corresponds to the Jacobian term $-J(\psi,\omega)$. We take $\nsur$, the dimension of the reduced-order state space, to be significantly smaller than the dimension of the space of the discretized streamfunction, $\nvar$. It is known from~\cite{Sandu_2021_MF-EnKF} that for this specific spatial discretization, utilizing $\nvar=8001$ points, $\nsur=40$ modes corresponds to $98.1\%$ of the total kinetic energy of the system, while $\nsur=80$ modes corresponds to $99.4\%$ of the total energy, meaning that this problem should be extremely well-suited to dimensional reduction.

The timespan for which we run the system is a ten day forecast, which for the given parameters is approximately $t \in [0, 0.109]$. The ROM was built in the interpolatory regime on the test timespan using 2000 evenly-spaced snapshots.  The model and the reduced-order model implementations for the experiments have been taken from~\cite{Sandu_2019_ODE-tests}.  As they are written in MATLAB, we use \texttt{ode45} with absolute and relative tolerances of $10^{-11}$ to accurately solve the modified fast ODEs \cref{eq:SM-MRI-GARK:ode,eq:SM-SPC-MRI-GARK:ode}.

\begin{figure}[ht!]
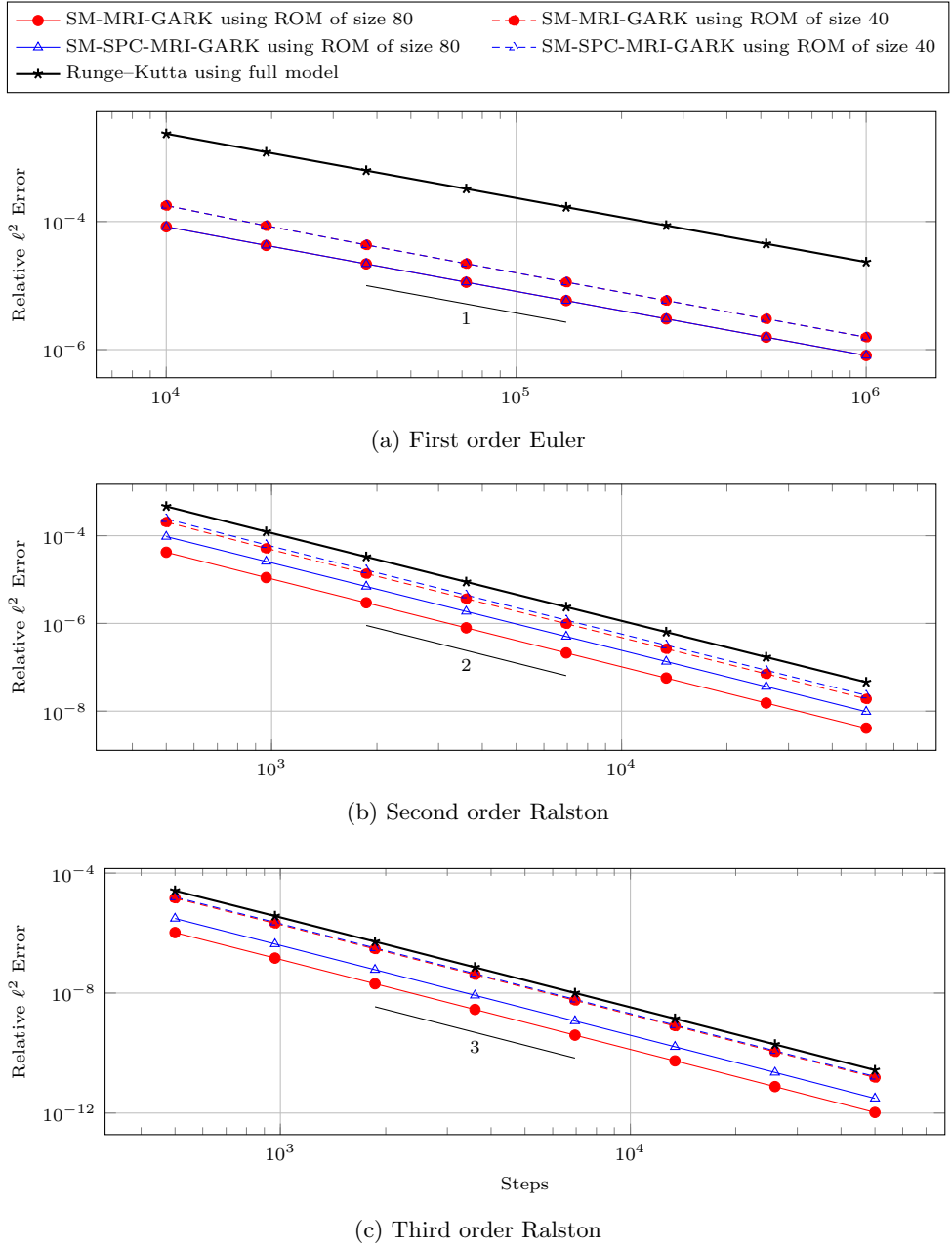

	\centering
	\begin{subfigure}{\linewidth}
		\centering
		\ref*{legend:qg_conv}
	\end{subfigure}

	\medskip
	
	\begin{subfigure}{\linewidth}
		\qgconvplot{Euler}{1e-5}
		\caption{First order Euler}
	\end{subfigure}

	\begin{subfigure}{\linewidth}
		\qgconvplot{Ralston2}{9e-7}
		\caption{Second order Ralston}
	\end{subfigure}

	\begin{subfigure}{\linewidth}
		\qgconvplot{Ralston3}{3.5e-9}*
		\caption{Third order Ralston}
	\end{subfigure}

	\caption{Convergence plots for the Quasi-Geostrophic equations \cref{eq:qge}.}
	\label{fig:qg-convergence}
\end{figure}

\Cref{fig:qg-convergence} shows the convergence of the QG equations with respect to a ROM.  We have not plotted the results of the Runge--Kutta methods using only the surrogate model as they produced errors orders of magnitude larger than the other methods.  For a fixed number of timesteps, \mrimethod{} and \spcmethod{} consistently have less error than Runge--Kutta methods using only the full model.  At order three, however, the more accurate ROM of size $\nsur=80$ is needed to achieve a substantial decrease in error.

\subsection{Lorenz '96 with a machine learning surrogate model}

The Lorenz '96 problem \cite{lorenz1996predictability} is given by
\begin{equation} \label{eq:Lorenz96}
    X'_k  = -X_{k-2} X_{k-1} + X_{k-1} X_{k+1} - X_k + F, \qquad k = 1, 2, \ldots, 40,
\end{equation}
 with periodic boundary conditions, and a forcing term $F = 8$. Due to its chaotic dynamics, this problem is prominently used in data assimilation and atmospheric research. Readers interested in current literature on machine learning surrogate models based on Lorenz '96 may refer to \cite{Rasp2020OnlineLorenz,Gagne2020Mar}. 
 
 In the following experiments, a neural network model with fully-connected layers was trained to approximate the right-hand side function in \cref{eq:Lorenz96}. Therefore the projections are $\Vmat = \Wmat* = \eye{d}$ and the  trained network acts as $\fsur $.  The network consists of  input and output layers and a hidden layer of dimension $80$. The input and hidden layers use the softplus activation function while the output layer does not have any activations.  Again, this experiment is implemented in MATLAB and uses \texttt{ode45} with tight tolerances to solve ODEs involving the surrogate model.
 
Similar to \cite{lorenz1996predictability}, the initial condition $X_{20}(0) = 8.008$ and $X_i(0) = 8$ for $i \ne 20$
%
%
is propagated for 2 units of time to expel transient effects and the developed solution is used as the true initial condition.  The training  data is taken from 5000 equally-spaced solutions of the full model within the interval $t \in [2,10]$ and, the convergence tests are performed over timespan $t \in [4, 8]$.

\begin{figure}[ht!]
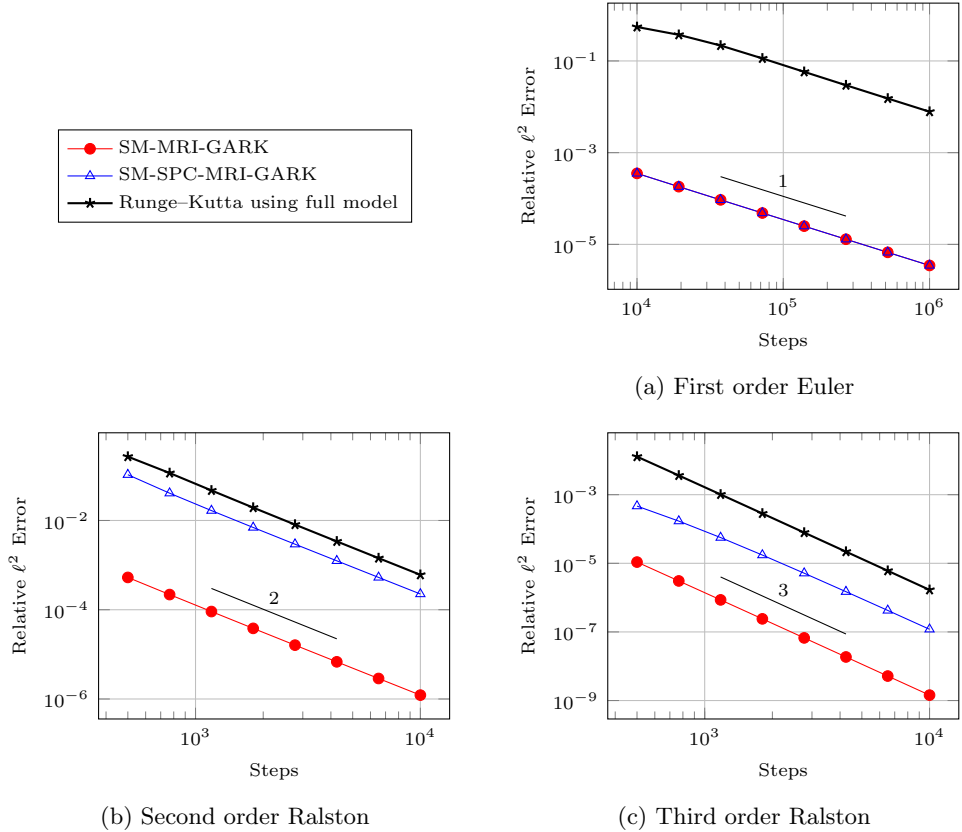

	\centering
	\begin{subfigure}[t]{0.48\linewidth}
		\centering
		\ref*{legend:ml_conv}
	\end{subfigure}
	\hfill
	\begin{subfigure}{0.48\linewidth}
		\mlconvplot{Euler}{3e-4}
		\caption{First order Euler}
	\end{subfigure}
	\hfill
	\begin{subfigure}{0.48\linewidth}
		\mlconvplot{Ralston2}{3e-4}
		\caption{Second order Ralston}
	\end{subfigure}
	\hfill
	\begin{subfigure}{0.48\linewidth}
		\mlconvplot{Ralston3}{4e-6}*
		\caption{Third order Ralston}
	\end{subfigure}
	\caption{Convergence plots for the Lorenz '96 problem \cref{eq:Lorenz96}.}
    \label{fig:convergence-ml-lorenz96}
\end{figure}

The convergence results are shown in \cref{fig:convergence-ml-lorenz96}. We note that merely integrating the surrogate model did not produce a stable solution in any of the experiments. On the other hand, when used together with the full model in \mrimethod{} and \spcmethod{} methods, the accuracy of the solution increases significantly compared to the full model solution.

\ifreport

\subsubsection{Stepsize adaptivity}


Outside of \cref{rem:adaptivity}, we have only considered fixed values of $H$ across an entire timespan of interest.  Some adaptivity has been introduced by using adaptive methods like \texttt{ode45} to solve modified fast ODEs.  For this section, we present $H$ adaptivity results for \mrimethod{} and \spcmethod{} methods applied to the Lorenz '96 problem.  We use the error controller from \cite[Section II.4]{Hairer1993} to select $H$ such that the local truncation error is within specified absolute and relative tolerances (\texttt{AbsTol} and \texttt{RelTol}, respectively).

For the first experiment, we test the second order \spcmethod{} and third order \mrimethod{} from \cref{app:method_coeffs} on the Lorenz '96 problem with its machine-learning-based surrogate.  \Cref{fig:adaptivity-ml-lorenz96} plots the adaptively chosen $H$ for each step.  We can see the value of $H$ varies slightly throughout the timespan and there are reasonably few rejected steps.

\begin{figure}[ht!]
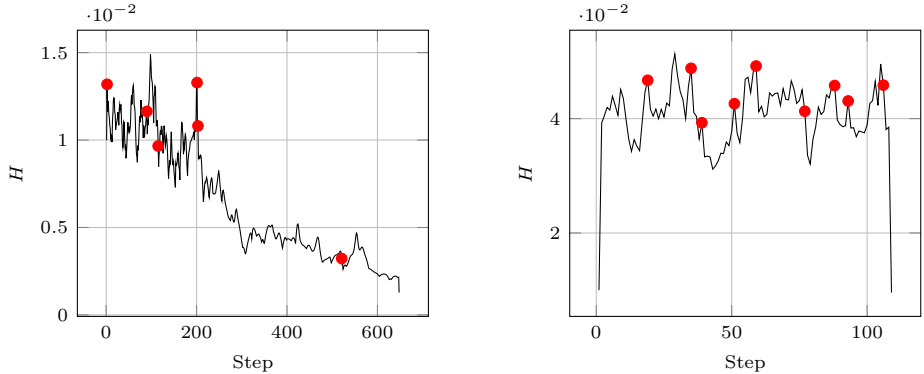

	\centering
	\begin{subfigure}[t]{0.48\linewidth}
		\mladaptplot{spc_Ralston2}
		\caption{Second order Ralston \spcmethod{}}
	\end{subfigure}
	\hfill
	\begin{subfigure}[t]{0.48\linewidth}
		\mladaptplot{gark_Ralston3}
		\caption{Third order Ralston \mrimethod{}}
	\end{subfigure}
	\caption{Adaptivity selected stepsize $H$ for each step taken to solve the Lorenz '96 problem \cref{eq:Lorenz96} with $\texttt{AbsTol}=\texttt{RelTol}=10^{-4}$.  Rejected steps shown with red markers.}
	\label{fig:adaptivity-ml-lorenz96}
\end{figure}

For the second experiment we test all methods in \cref{app:method_coeffs} using $\texttt{AbsTol} = \texttt{RelTol} = 10^{-3}, 10^{-4}, \dots, 10^{-10}$.  The (global) error at the end of the timespan is recorded for each tolerance and plotted in \cref{fig:tol-ml-lorenz96}.  For reference, we also include traditional, adaptive Runge--Kutta methods using only the full model.  As expected, cutting the local error tolerance by a factor of ten causes the global error to decrease by a factor of approximately ten.
\begin{figure}[ht!]
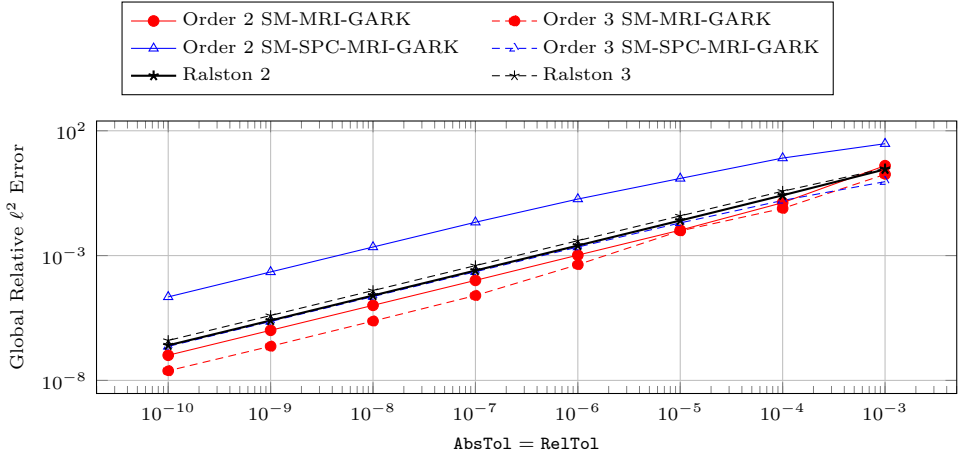

	\centering
	\begin{subfigure}{\linewidth}
		\centering
		\ref*{legend:ml_tol}
	\end{subfigure}
	
	\medskip
	
	\mltolplot{}
	
	\caption{Global error versus stepsize controller tolerances for the Lorenz '96 problem \cref{eq:Lorenz96}.}
	\label{fig:tol-ml-lorenz96}
\end{figure}

\fi

\subsection{Brusselator PDE}
\label{subsec:brusselator}

The test problem BRUS from \cite[Section II.10]{Hairer1993} is based on the Brusselator reaction-diffusion PDE 
\begin{equation} \label{eq:brus}
	\begin{split}
		\pdv{u}{t} &= \alpha \laplacian{u} + 1 + u^2 v - 4.4 u, \\
		\pdv{v}{t} &= \alpha \laplacian{v} + 3.4 u - u^2 v,
	\end{split}
\end{equation}
where $\alpha = 0.002$.  The spatial domain is the unit square $0 \le x, y \le 1$ with zero, Neumann boundary conditions for both $u$ and $v$. Using the method of lines, this is discretized with second order central finite difference on a uniform $P \times P$ grid.  Starting at the initial conditions
\begin{equation*}
	u(t=0, x, y) = 0.5 + y,
	\qquad
	v(t=0, x, y) = 1 + 5 x,
\end{equation*}
we seek the solution at $t = 7.5$.

In our experiments, $P = 257$, and thus, $N = 2 P^2 = 132098$.  For the surrogate model, we take the approach described in \cref{subsec:multimesh} where we use the same finite difference discretization of \cref{eq:brus} but on a coarser mesh.  In particular, coarse meshes of size $P = 129$ and $P = 65$ are used, which nest inside the fine mesh.  The modified fast ODEs \cref{eq:SM-MRI-GARK:ode,eq:SM-SPC-MRI-GARK:ode} are solved with \textit{one} step of a Runge--Kutta method one or two orders higher than the base method.  We found that this uses the fewest evaluations of $\fsur$ while keeping the ODE solution errors negligible.

\begin{figure}[ht!]
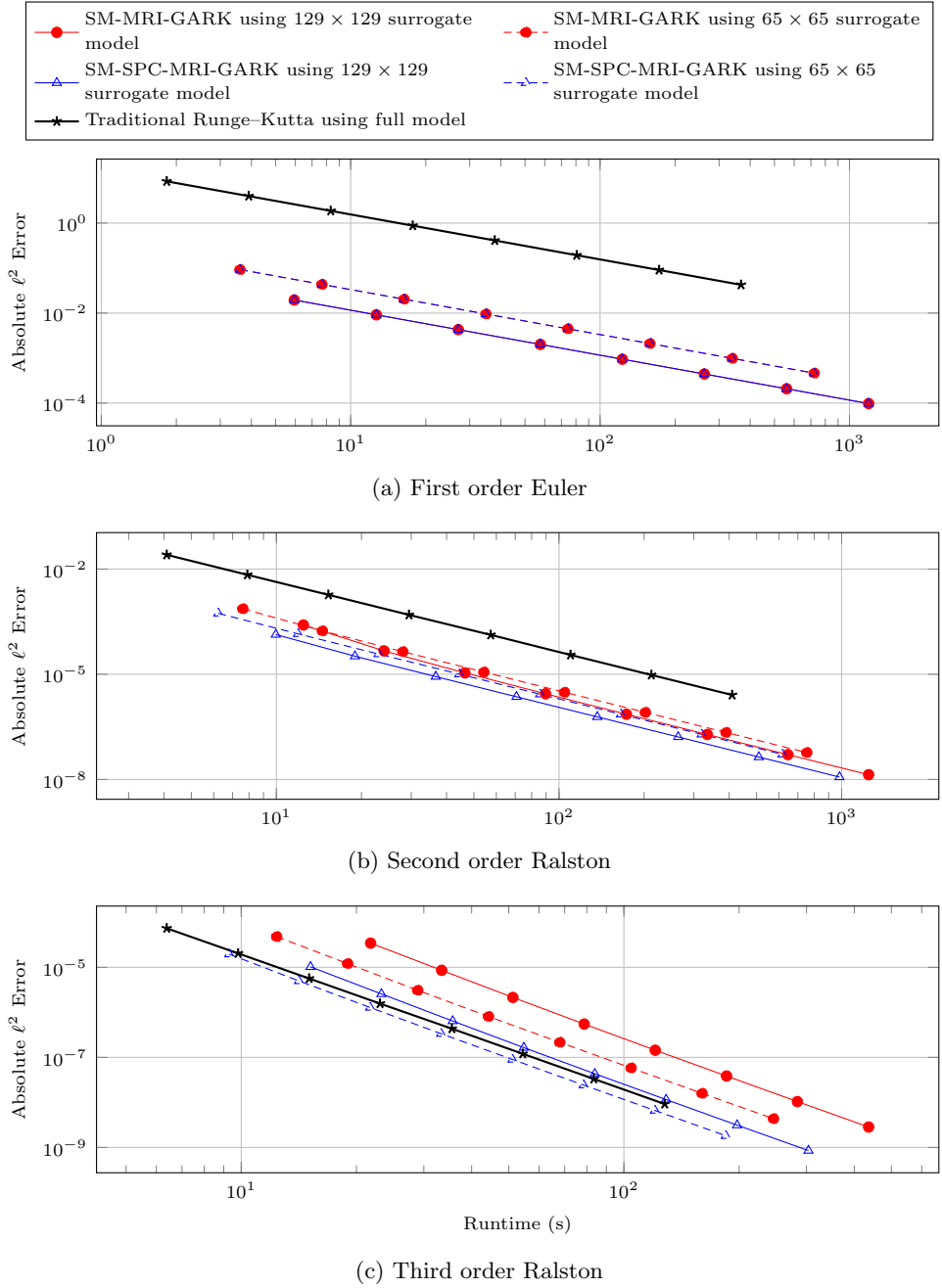

	\centering
	\begin{subfigure}{\linewidth}
		\centering
		\ref*{legend:brus_perf}
	\end{subfigure}

	\medskip

	\begin{subfigure}{\linewidth}
		\brusperfplot{1}
		\caption{First order Euler}
	\end{subfigure}

	\begin{subfigure}{\linewidth}
		\brusperfplot{2}
		\caption{Second order Ralston}
	\end{subfigure}

	\begin{subfigure}{\linewidth}
		\brusperfplot{3}*
		\caption{Third order Ralston}
	\end{subfigure}
	\caption{Work precision diagrams for BRUS \cref{eq:brus}.}
	\label{fig:brus}
\end{figure}

Comparing the performance of the integrators is our primary goal with the Brusselator problem, and we implemented the tests in C.  For each integrator and surrogate model size, the runtime and error was recorded for a range of eight stepsizes.  The results are plotted in \cref{fig:brus}.  At orders one and two, \mrimethod{} and \spcmethod{} show clear speedups over Runge--Kutta solutions using only the full model.  In addition, the performance increases as the surrogate model mesh becomes finer.  Unfortunately, the results are reversed at order three.

Profiling of the code helped to identify why \mrimethod{} and \spcmethod{} performed poorly.  \Cref{tab:brus_stats} provides a breakdown of the relative costs of four, primary operations.  When a $127 \times 127$ surrogate model is used, about one third of the runtime is devoted to solving modified fast ODEs, and an evaluation of $\fsur$ is about 25\% as expensive as $\fun$.  The linear operators $\Vmat$ and $\Wmat*$ have efficient, matrix-free implementations, but are still 75\% as expensive as evaluating $\fun$.  These factors caused one step of \mrimethod{} and \spcmethod{} to take approximately twice as long as a traditional Runge--Kutta step.  At order three, the reduction in error is not enough to overcome the overheads.

\begin{table}
	\centering
	\begin{tabular}{r||c|c|c|c|c|c|c|c|c|c|c|c}
		Method & \multicolumn{6}{c|}{\mrimethod{}} & \multicolumn{6}{c}{\spcmethod{}} \\ \hline
		Surrogate model &
		\multicolumn{3}{c|}{$129 \times 129$} & \multicolumn{3}{c|}{$65 \times 65$} & \multicolumn{3}{c|}{$129 \times 129$} & \multicolumn{3}{c}{$65 \times 65$} \\ \hline
		Order & 1 & 2 & 3 & 1 & 2 & 3 & 1 & 2 & 3 & 1 & 2 & 3 \\ \hline \hline
		$\fun$ evals & 27 & 26 & 22 & 44 & 42 & 39 & 27 & 32 & 32 & 44 & 52 & 52 \\
		Modified fast ODEs & 30 & 31 & 38 & 12 & 12 & 16 & 30 & 29 & 33 & 12 & 10 & 12 \\
		$\Vmat$ and $\Wmat*$ evals & 25 & 23 & 20 & 27 & 25 & 23 & 24 & 16 & 11 & 27 & 16 & 11 \\
		Other & 18 & 20 & 20 & 18 & 20 & 22 & 18 & 23 & 25 & 18 & 22 & 25
	\end{tabular}
	\caption{Percentage of total runtime devoted to four suboperations for \mrimethod{} and \spcmethod{} methods applied to BRUS \cref{eq:brus}.  ``Other'' includes, for example, the time spent allocating memory and computing linear combinations of stages.  Note that each column sums to $100\%$.}
	\label{tab:brus_stats}
\end{table}

\subsection{Advection PDE}

Finally, we consider the following 2D advection problem with zero, Dirichlet boundary conditions:
\begin{alignat*}{2}
	\pdv{u}{t} + a \cdot \nabla{u} &= 0, & \qquad & \text{on } \Omega = [0, 1] \times [0, 1], \\
	u(t,x,y) &= 0, & \qquad & \text{on } \partial \Omega.
\end{alignat*}
The velocity field corresponds to the Molenkamp-Crowley problem:
\begin{equation*}
	a(x, y) = \left[
		2 \pi \left(y - \sfrac{1}{2}\right),
		-2 \pi \left(x - \sfrac{1}{2}\right)
	\right].
\end{equation*}
We start with the initial condition
\begin{equation*}
	u(t = 0, x, y) = \exp( -100 \left( (x - 0.35)^2 + (y - 0.35)^2 \right))
\end{equation*}
and stop at $t = 2$.  This allows the initial profile to rotate clockwise twice about the center of the domain.

For the method of lines discretization in space, we represent $\Omega$ with a $P \times P$ uniform, triangular mesh and apply a nodal discontinuous Galerkin (DG) method.  The order in space is chosen to match the order in time. This yields the linear ODE
\begin{equation} \label{eq:advection_ode}
	\y' = \mathbf{M}^{-1}\, \mathbf{K}\, \y,
\end{equation}
where $\mathbf{M}$ and $\mathbf{K}$ are mass and advection matrices, respectively.  We use the C++ library MFEM \cite{mfem} for this DG discretization, and in particular, base it on the serial version of example 9 provided in MFEM version 4.1.

Again, we try the multimesh approach by using a mesh size of $P=100$ for $\fun$ and a mesh size of $P=50$ for $\fsur$.  We use the same technique discussed in \cref{subsec:brusselator} for integrating modified fast ODEs involving $\fsur$.  The size of the ODE \cref{eq:advection_ode} for each experiment is listed in \cref{tab:advection_size}.

\begin{table}[ht!]
	\centering
	\begin{tabular}{c|c|c}
		\shortstack{Order of time and\\space discretization} & \shortstack{Dimension of\\full model $\nvar$} & \shortstack{Dimension of\\surrogate model $\nsur$} \\ \hline
		1 & $6 \times 10^4$ &  $1.5 \times 10^4$ \\
		2 & $1.2 \times 10^5$ &  $3 \times 10^4$ \\
		3 &  $2 \times 10^5$ &  $5 \times 10^4$ \\
	\end{tabular}
	\caption{Dimensions of the full and surrogate models used in the advection experiment.}
	\label{tab:advection_size}
\end{table}

In contrast to the Brusselator problem \cref{eq:brus}, the advection problem is linear and hyperbolic.  Moreover, profiling reveals an evaluation of the RHS of \cref{eq:advection_ode} is much more expensive than an interpolation between meshes due, in part, to the linear solve with the mass matrix.  Based on the discussion in the previous subsection, we may expect better performance results at order three, and indeed, that is the case.

\begin{figure}[ht!]
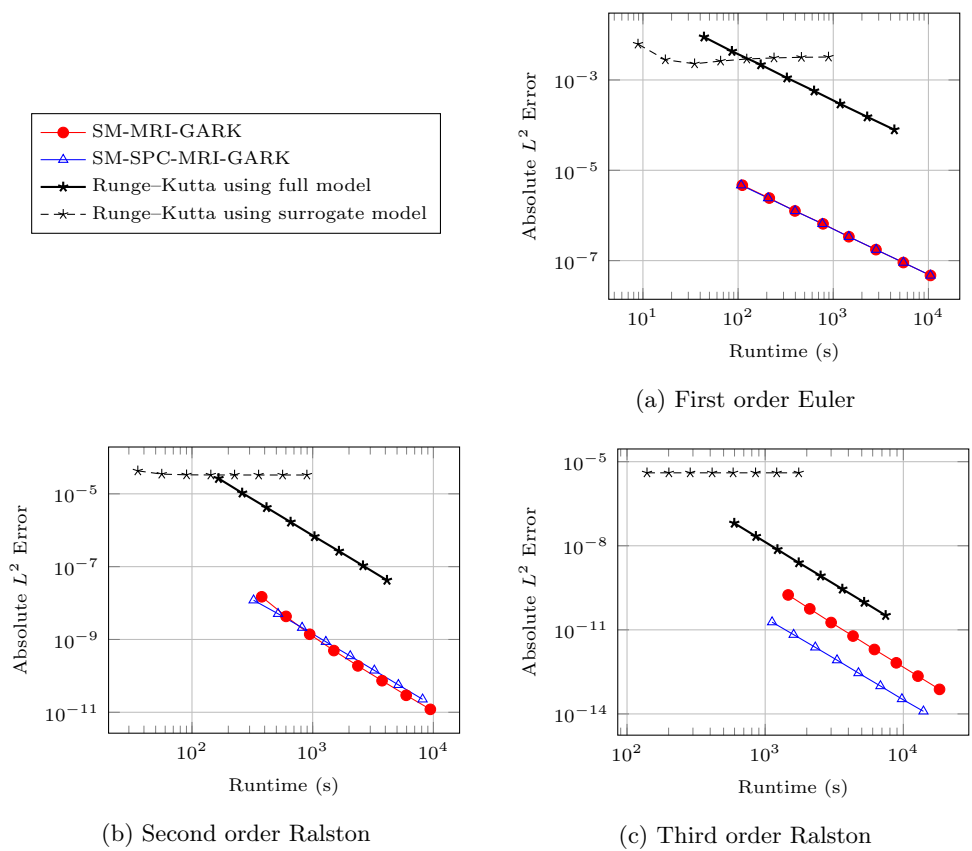

	\centering
	\begin{subfigure}[t]{0.48\linewidth}
		\centering
		\ref*{legend:adv_perf}
	\end{subfigure}
	\hfill
	\begin{subfigure}{0.48\linewidth}
		\advperfplot{1}
		\caption{First order Euler}
	\end{subfigure}
	\hfill
	\begin{subfigure}{0.48\linewidth}
		\advperfplot{2}
		\caption{Second order Ralston}
	\end{subfigure}
	\hfill
	\begin{subfigure}{0.48\linewidth}
		\advperfplot{3}*
		\caption{Third order Ralston}
	\end{subfigure}
    \caption{Work precision diagrams for advection problem \cref{eq:advection_ode}.}
    \label{fig:advection}
\end{figure}

\Cref{fig:advection} shows the error and timing for the integrators over a range of eight stepsizes.  All \mrimethod{} and \spcmethod{} methods outperform Runge--Kutta methods with speedups ranging from approximately 3 to 725.  The error for the Runge--Kutta methods using the surrogate model remains nearly constant, which indicates that spatial errors dominate temporal errors.

\section{Conclusions}
\label{sec:conclusion}

This work develops a new methodology to accelerate the time integration of large ODEs using surrogate models. 
Specifically, surrogate information is incorporated into the numerical solution of the original ODEs using multirate methods.  We derive and analyze two implementations of this new time-stepping technique: \mrimethod{} and \spcmethod{}.  Both combine continuous evolution of the surrogate model with discrete Runge--Kutta steps of the full ODE.  There are numerous ways to generate surrogate models that offer flexible trade-offs among accuracy, computational cost, and stability.  The new methods are designed such that the overall order of accuracy is independent of the surrogate model. The more accurate the surrogate models are, the smaller the local error constants are, which leads to a smaller global error.

Numerical experiments reveal that, at low orders of accuracy, it is possible to achieve orders of magnitude speedups over traditional Runge--Kutta methods.  As the order increases, overheads associated with evaluations of $\Vmat$, $\Wmat*$, and the surrogate model are penalized more when measuring efficiency, and speedups tend to decrease.  There is not a clear winner between \mrimethod{} and \spcmethod{}; their relative performance appears to depend heavily on the properties of the full and surrogate models.

The methods considered in this paper are based on explicit multirate infinitesimal methods, although the surrogate integration can be done implicitly.  Even so, the explicit treatment of the surrogate model error can lead to stepsize restrictions.  Instability can occur if a surrogate model does not capture the stiffness of the full model.  The authors plan to explore implicit and implicit-explicit methods \cite{Sandu_2014_IMEX-GLM,Sandu_2014_IMEX_GLM-Extrap,Sandu_2015_IMEX-TSRK,Sandu_2016_highOrder-IMEX-GLM} for stiff problems and differential-algebraic equations.  We note that stiffness can appear in the full model, surrogate model, and even the surrogate model error, and each scenario brings different considerations.  

\appendix
\section{Explicit MRI-GARK and SPC-MRI-GARK methods of orders two and three}
\label{app:method_coeffs}

\begin{figure}[ht!]
	\centering
	\begin{subfigure}{0.48\linewidth}
		\includegraphics[width=\linewidth]{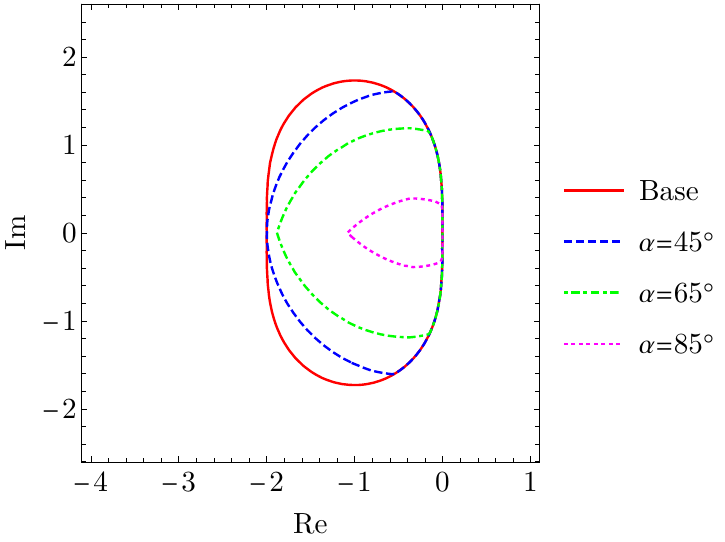}
		\caption{Second order MRI-GARK}
	\end{subfigure}
	\hfill
	\begin{subfigure}{0.48\linewidth}
		\includegraphics[width=\linewidth]{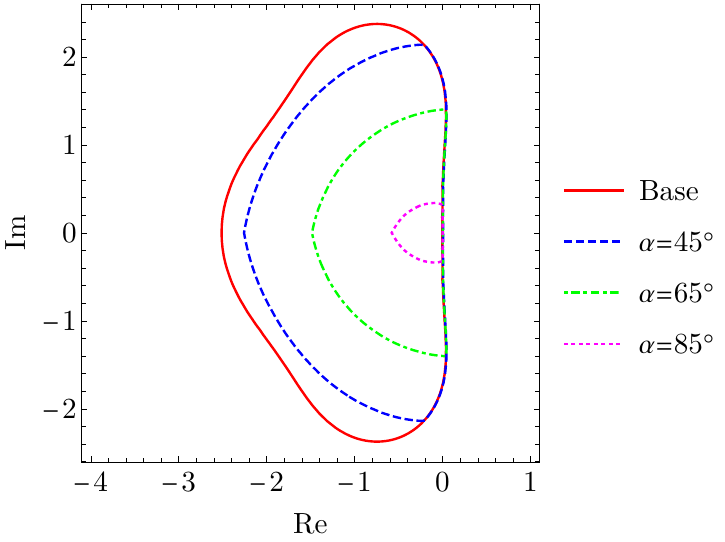}
		\caption{Third order MRI-GARK}
	\end{subfigure}
	\\ \medskip
	\begin{subfigure}{0.48\linewidth}
		\includegraphics[width=\linewidth]{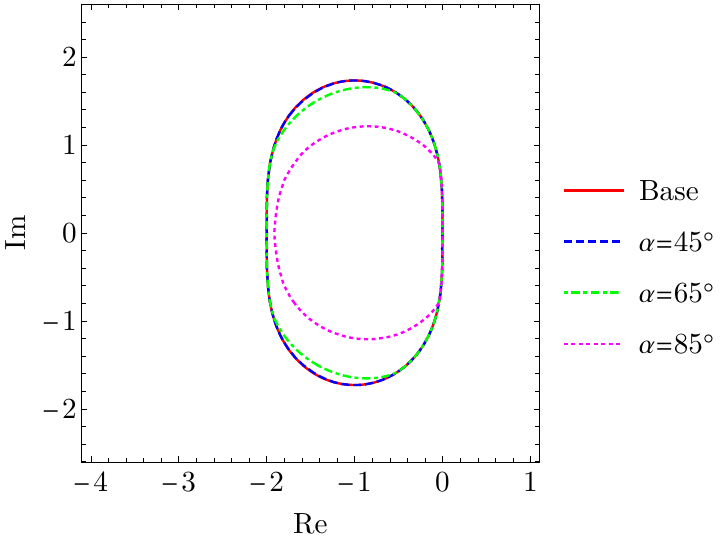}
		\caption{Second order SPC-MRI-GARK}
	\end{subfigure}
	\hfill
	\begin{subfigure}{0.48\linewidth}
		\includegraphics[width=\linewidth]{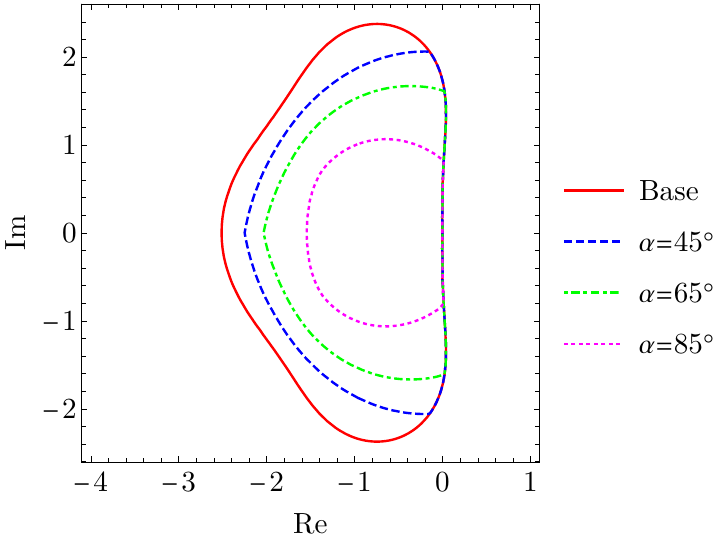}
		\caption{Third order SPC-MRI-GARK}
	\end{subfigure}
	\caption{Stability regions for new methods in \cref{tab:coeffs_order} including the base Runge--Kutta stability region and $\scalarstab[\infty]$ for $\alpha = \ang{45}, \ang{65}, \ang{85}$.}
	\label{fig:stability}
\end{figure}

\Cref{tab:coeffs_order} provides coefficients for new MRI-GARK and SPC-MRI-GARK schemes based on Ralston's optimal second and third order Runge--Kutta methods \cite{ralston1962runge}.  \Cref{fig:stability} plots their scalar slow stability regions \cite[Definition 4.1]{Sandu_2019_MRI-GARK}
\begin{equation*}
	\scalarstab[\infty] = \big\{ z\S \in \Cplx \, \big| \, \big|R \big( z\F, z\S \big)\big| \le 1, \, \forall z\F \in \Cplx{-} : \big|\text{arg}{\big( z\F \big)} - \pi \big| \le \alpha \big\},
\end{equation*}
where $R(z\F, z\S)$ is the scalar linear stability function.  All methods in \cref{tab:coeffs_order} satisfy $\lim_{z\F\to-\infty} R(z\F, z\S) = 0$, so they are suitable for problems where the fast dynamics are stiff, but the slow dynamics are nonstiff.  For readers interested in using different base methods, we have provided a Mathematica notebook in the supplementary materials with general, parameterized families of methods.

\begin{table}[ht!]
	\centering
	\begin{tabular}{c||c||c||c}
		Order & Base Method & MRI-GARK $\Gamma(t), \widehat{\gamma}(t)$ & \shortstack{SPC-MRI-GARK\\$\gamma(t), \widehat{\gamma}(t)$} \\ \hline \hline
		2 &
		$\begin{butchertableau}{c|cc}
		0 & 0 & 0 \\
		\frac{2}{3} & \frac{2}{3} & 0 \\ \hline
		& \frac{1}{4} & \frac{3}{4}
		\end{butchertableau}$
		&
		$\begin{smatrix}
			\frac{2}{3} & 0 \\
			-\frac{5}{12} & \frac{3}{4}
		\end{smatrix}, \begin{smatrix}
			\frac{1}{3} \\ 0
		\end{smatrix}$
		&
		$\begin{smatrix}
			-\frac{1}{2} + \frac{3 t}{2} \\
			\frac{3}{2} - \frac{3 t}{2}
		\end{smatrix}, \begin{smatrix}
			1 \\ 0
		\end{smatrix}$ \\ \hline \hline
		3 &
		$\begin{butchertableau}{c|ccc}
			0 & 0 & 0 & 0 \\
			\frac{1}{2} & \frac{1}{2} & 0 & 0 \\
			\frac{3}{4} & 0 & \frac{3}{4} & 0 \\ \hline
			& \frac{2}{9} & \frac{1}{3} & \frac{4}{9}
		\end{butchertableau}$
		&
		$\substack{\begin{smatrix}
			\frac{1}{2} & 0 & 0 \\
			-\frac{11}{4} + \frac{9 t}{2} & 3-\frac{9 t}{2} & 0 \\
			\frac{47}{36}-\frac{13 t}{6} & -\frac{1}{6}-\frac{t}{2} & -\frac{8}{9}+\frac{8 t}{3}
		\end{smatrix}, \\ \begin{smatrix}
			\frac{1}{40} \\ \frac{7}{40} \\ \frac{1}{20}
		\end{smatrix}}$
		&
		$\substack{\begin{smatrix}
			1-\frac{2 t}{3}-\frac{4 t^2}{3} \\
			-2 t + 4 t^2 \\
			\frac{8 t}{3}-\frac{8 t^2}{3}
		\end{smatrix}, \\ \begin{smatrix}
			-\frac{7}{8} + \frac{9 t}{5} \\ \frac{71}{40}-\frac{17 t}{10} \\ \frac{1}{10}-\frac{t}{10}
		\end{smatrix}}$
	\end{tabular}
	\caption{Second and third order MRI-GARK and SPC-MRI-GARK coefficients.}
	\label{tab:coeffs_order}
\end{table}

\bibliographystyle{siamplain}
\bibliography{Bib/main,Bib/ode_multirate,Bib/sandu}

\end{document}